%------------------------------------------------------------------------------

\magnification 1200
\input amstex.tex
\input amsppt.sty

\voffset=-36pt

\author Nazarov and Tarasov \endauthor
\title Tensor Products of Yangian Modules \endtitle

%------------------------------------------------------------------------------

\expandafter\ifx\csname maxim.def\endcsname\relax \else\endinput\fi
\expandafter\edef\csname maxim.def\endcsname{%
 \catcode`\noexpand\@=\the\catcode`\@\space}
\catcode`\@=11

\mathsurround 1.6pt
\font\Bbf=cmbx12 

\def\hcor#1{\advance\hoffset by #1}
\def\vcor#1{\advance\voffset by #1}
\let\bls\baselineskip  \let\ignore\ignorespaces
\def\vsk#1>{\vskip#1\bls} \let\adv\advance 
\def\vv#1>{\vadjust{\vsk#1>}\ignore} \def\vvv#1>{\vadjust{\vskip#1}\ignore}
\def\vvn#1>{\vadjust{\nobreak\vsk#1>\nobreak}\ignore}
\def\vvvn#1>{\vadjust{\nobreak\vskip#1\nobreak}\ignore}
\def\emph#1{{\it #1\/}}

\let\vp\vphantom  %\let\^\negthickspace
 \let\nt\noindent 
\def\nn#1>{\noalign{\vskip #1pt}} \def\NN#1>{\openup#1pt}
 
\let\Sum\sum \def\sum{\Sum\limits} 
\let\Prod\prod \def\prod{\Prod\limits} \let\Int\int \def\int{\Int\limits}

\let\=\m@th \def\&{.\kern.1em} \def\>{\!\;} \def\:{\!\!\;}

\ifx\plainfootnote\undefined \let\plainfootnote\footnote \fi
\expandafter\ifx\csname amsppt.sty\endcsname\relax
 
\else \fi

\newbox\s@ctb@x
\def\s@ct#1 #2\par{\removelastskip\vsk>
 \vtop{\bf\setbox\s@ctb@x\hbox{#1} \parindent\wd\s@ctb@x
 \ifdim\parindent>0pt\adv\parindent.5em\fi\item{#1}#2\strut}%
 \nointerlineskip\nobreak\vtop{\strut}\nobreak\vsk-.4>\nobreak}

\newbox\t@stb@x
\def\gadv{\global\advance} \def\gad#1{\gadv#1 1} 
\def\l@b@l#1#2{\def\n@@{\csname #2no\endcsname}%
 \if *#1\gad\n@@ \expandafter\xdef\csname @#1@#2@\endcsname{\the\Sno.\the\n@@}%
 \else\expandafter\ifx\csname @#1@#2@\endcsname\relax\gad\n@@
 \expandafter\xdef\csname @#1@#2@\endcsname{\the\Sno.\the\n@@}\fi\fi}
\def\l@bel#1#2{\l@b@l{#1}{#2}\?#1@#2?}
\def\?#1?{\csname @#1@\endcsname}
\def\[#1]{\def\n@xt@{\ifx\t@st *\def\n@xt####1{{\setbox\t@stb@x\hbox{\?#1@F?}%
 \ifnum\wd\t@stb@x=0 {\bf???}\else\?#1@F?\fi}}\else
 \def\n@xt{{\setbox\t@stb@x\hbox{\?#1@L?}\ifnum\wd\t@stb@x=0 {\bf???}\else
 \?#1@L?\fi}}\fi\n@xt}\futurelet\t@st\n@xt@}
\def\(#1){{\rm\setbox\t@stb@x\hbox{\?#1@F?}\ifnum\wd\t@stb@x=0 ({\bf???})\else
 (\?#1@F?)\fi}}
\def\dff{\expandafter\d@f} \def\d@f{\expandafter\def}
\def\edff{\expandafter\ed@f} \def\ed@f{\expandafter\edef}

\newcount\Sno \newcount\Lno \newcount\Fno
\def\Section#1{\gad\Sno\Fno=0\Lno=0\s@ct{\the\Sno.} {#1}\par} \let\Sect\Section
\def\section#1{\gad\Sno\Fno=0\Lno=0\s@ct{} {#1}\par} \let\sect\section
\def\l@F#1{\l@bel{#1}F} \def\<#1>{\l@b@l{#1}F} \def\l@L#1{\l@bel{#1}L}
\def\Tag#1{\tag\l@F{#1}} \def\Tagg#1{\tag"\llap{\rm(\l@F{#1})}"}
\def\Th#1{Theorem \l@L{#1}} \def\Lm#1{Lemma \l@L{#1}}
\def\Prop#1{Proposition \l@L{#1}}
\def\Cr#1{Corollary \l@L{#1}} \def\Cj#1{Conjecture \l@L{#1}}
 
\def\Proof#1.{\demo{\it Proof #1}}

 \def\setparindent{\edef\Parindent{\the\parindent}}
\def\Appendix{\Sno=64\let\p@r@\z@ %\parindent
\def\Section##1{\gad\Sno\Fno=0\Lno=0 \s@ct{} \hskip\p@r@ Appendix \char\the\Sno
 \if *##1\relax\else {.\enspace##1}\fi\par} \let\Sect\Section
\def\section##1{\gad\Sno\Fno=0\Lno=0 \s@ct{} \hskip\p@r@ Appendix%
 \if *##1\relax\else {.\enspace##1}\fi\par} \let\sect\section
\def\l@b@l##1##2{\def\n@@{\csname ##2no\endcsname}%
 \if *##1\gad\n@@
\expandafter\xdef\csname @##1@##2@\endcsname{\char\the\Sno.\the\n@@}%
\else\expandafter\ifx\csname @##1@##2@\endcsname\relax\gad\n@@
 \expandafter\xdef\csname @##1@##2@\endcsname{\char\the\Sno.\the\n@@}\fi\fi}}

\let\logo@\relax
\let\m@k@h@@d\makeheadline \let\m@k@f@@t\makefootline
\def\makeheadline{\ifnum\pageno=1\headline={\hfil}\fi\m@k@h@@d}
\def\makefootline{\ifnum\pageno=1\footline={\hfil}\fi\m@k@f@@t}

%------------------------------------------------------------------------------

\def\al{\alpha}
\def\AN{\operatorname{A}(\glN)}

\def\be{\beta}
\def\bi{{\bold i}}
\def\bj{{\bold j}}
\def\bk{{\bold k}}
\def\bx{{\boxed{\phantom{\square}}\kern-.4pt}}

\def\CC{{\Bbb C}}

\def\de{\delta}
\def\De{\Delta}
\def\deg{{\operatorname{deg}\ts}}
\def\Den{{\De^{\hskip-1pt(n)}}}
\def\Dep{{\Delta^{\hskip-1pt\prime}}}
\def\Der{{\De^{\hskip-1pt(r)}}}

\def\End{{\operatorname{End}}}
\def\EndCN{\End(\CC^N)}
\def\enddemos{{$\hfill\square$\enddemo}}

\def\ga{\gamma}
\def\Gal{{\Gamma_{\hskip-1pt\al\be}}}
\def\ge{\geqslant}
\def\gl{{\frak{gl}}}
\def\glM{{\gl_M}}
\def\glMN{{\gl_{M+N}}}
\def\glN{{\gl_N}}

\def\id{{\operatorname{id}}}
\def\im{{\operatorname{im}\ts}}
\def\io{\iota}

\def\ka{\kappa}
\def\ker{{\operatorname{ker}\ts}}

\def\la{\lambda}
\def\La{\Lambda}
\def\lcd{{\ts,\hskip.95pt\ldots\ts,\ts}}
\def\le{\leqslant}
\def\lm{{\la,\mu}}
\def\lsm{{\la/\mu}}

\def\mi{{\raise.5pt\hbox{-}}}
\def\mv{{\kern55pt}}
\def\mw{{\kern84pt}}

\def\o{{\raise.24ex\hbox{$\sssize\kern.1em\circ\kern.1em$}}}
\def\of{{\om_f}}
\def\om{\omega}
\def\Om{\Omega}
\def\ox{\otimes}

\def\ph{\varphi}
\def\ps{\psi}

\def\R{{\Cal R}}
\def\ru{^{(r)}}

\def\S{{\Cal S}}
\def\si{\sigma}
\def\su{^{(s)}}

\def\th{{\tau_h}}
\def\tm{{\hskip-1pt-}}
\def\tp{{\hskip-1pt+}}
\def\ts{{\hskip1pt}}

\def\x{\times}
\def\X{{\Cal X}}

\def\UM{{\operatorname{U}(\glM)}}
\def\UMN{{\operatorname{U}(\glMN)}}
\def\UN{{\operatorname{U}(\glN)}}

\def\YMN{\operatorname{Y}(\glMN)}
\def\YN{\operatorname{Y}(\glN)}

\def\ze{\zeta}
\def\zes{\ze^\ast}
\def\ZZ{{\Bbb Z}}

%------------------------------------------------------------------------------

\line{\Bbf\hfill On Irreducibility of Tensor Products of Yangian Modules}
\bigskip\bigskip
\line{\Bbf\hfill Maxim Nazarov and Vitaly Tarasov}
\bigskip\bigskip

%------------------------------------------------------------------------------

\line{\bf Introduction\hfill}
\medskip\smallskip\nt
In the present article we continue our study [NT] of the finite-dimensional
modules over the Yangian $\YN$ of the general linear Lie algebra $\glN$.
\!The \text{Yangian $\!\YN\!$} is a canonical deformation of the universal
enveloping algebra $\operatorname{U}(\glN[u])$
in the class of Hopf algebras [D1]. Definition of the algebra $\YN$
in terms of an infinite family of generators $T_{ij}^{(s)}$
with $s=1,2\ts,\ldots$ and $i\ts,j=1\lcd N$ is given by (1.1),\,(1.2).
Comultiplication $\De:\YN\to\YN\ox\YN$ is defined by (1.1),\,(1.16).

The algebra $\YN$ admits an alternative definition
in terms of the ascending chain
$
\operatorname{U}(\frak{gl}_1)\subset\operatorname{U}(\frak{gl}_2)\subset\ldots
$
of classical universal enveloping algebras [O]. For any non-negative
integer $M$ consider the commutant in $\UMN$ of the subalgebra $\UM$.
This commutant is generated by the centre of the subalgebra $\UM$ and
a homomorphic image of the Yangian $\YN$, see Proposition 1.1.
The intersection of the kernels in $\YN$ of all these homomorphisms
for $M=0,1,2\ts,\ldots$ is zero.

For any dominant integral weights $\la$ and $\mu$
of the Lie algebras $\glMN$ and $\glM$ consider the subspace $V_{\lm}$
in the irreducible {\=$\glMN$}-module $V_\la$ of highest weight $\la$ formed
by all singular vectors with respect to $\glM$ of weight $\mu$.
The algebra $\YN$ acts in $V_{\lm}$ irreducibly through the above homomorphism.
Further, for any complex number $h$ there is an automorphism $\th$ of the
algebra $\YN$ defined in terms of the generating series (1.1) by the assignment
$T_{ij}(u)\mapsto T_{ij}(u+h)$. By pulling back the {\=$\YN$}-module
$V_{\lm}$ through this automorphism we obtain an irreducible {\=$\YN$}-module,
which we denote by $V_{\lm}(h)$ and call elementary.

Any formal series $f(u)\in 1+u^{-1}\ts\CC\ts[[u^{-1}]]$ also defines
an automorphism $\of$ of the algebra $\YN$ by
$T_{ij}(u)\mapsto f(u)\cdot T_{ij}(u)$.
Further, there is a canonical chain of algebras
$
\operatorname{Y}(\frak{gl}_1)\subset\ldots\subset\operatorname{Y}(\frak{gl}_N).
$
The elementary modules are distinguished amongst all irreducible
finite-dimensional {\=$\YN$}-modules $W$ by the following fact.
Take the commutative subalgebra $\AN$ in $\YN$ generated by the centres
of all algebras in the latter chain. Then the action of this subalgebra in $W$
is semi-simple, if and only if $W$ is obtained by pulling back through some
automorphism $\of$ from a tensor product of elementary {\=$\YN$}-modules
$$
V=V_{\la^{(1)},\,\mu^{(1)}}(h^{(1)})
\ox\ldots\ox
V_{\la^{(n)},\,\mu^{(n)}}(h^{(n)})
$$
where $h\ru-h\su\notin\ZZ$ whenever $1\le r<s\le n$.
This fact was conjectured in [C2] and proved in [NT]\ts.
It was subsequently applied in [KKN] and [TU] to the analysis of
integrable lattice models.
If $h\ru-h\su\notin\ZZ$ for all $r<s$ then
{\=$\YN$}-module $V$ is irreducible and the action of $\AN$ in $V$
has a simple spectrum.
Here we study the tensor product $V$ when
the differences $h\ru-h\su$ are any complex numbers.

Theorem 3.3 gives sufficient conditions for irreducibility of
{\=$\YN$}-module $V\!$.
In general, these conditions are not necessary for the irreducibility of $V\!$.
However, in the particular case when $\la^{(1)}\lcd\la^{(n)}$ are any
fundamental weights of $\glN$ whilst $\mu^{(1)}\lcd\mu^{(n)}$ are empty,
the conditions of Theorem 3.3 are necessary for the irreducibility of $V$
as well. This follows from Theorem 2.3. In this particular case the
irreducibility criterion of $V$ was obtained in [AK] by using the technique
of crystal bases. Our approach is different. We use an eigenbasis for the
action of $\AN$ in each tensor factor $V_{\la\su,\,\mu\su}(h\su)$ of $V\!$.
This eigenbasis is a part of the Gelfand-Zetlin basis [GZ] in the
{\=$\frak{gl}_{M\su+N}$}-module $V_{\la\su}$
corresponding to the chain of Lie algebras
$
\frak{gl}_1\subset\frak{gl}_2\subset\ldots\subset\ts\frak{gl}_{M\su+N}\,\ts.
$
Note that
a vector $\ze\su\in V\su$ called singular is contained in this eigenbasis.

Recall the notion of a universal {\=$R$}-matrix [D1] for the Hopf algebra
$\YN$. Let $\Dep$ be the composition of the comultiplication $\De$ with
permutation of the tensor factors in $\YN\ox\YN$. There exists a formal power
series $\R(z)$ in $z^{-1}$ with coefficients from $\YN\ox\YN$ and the leading
term $1$ such that
$$
\R(z)\cdot\id\ox\tau_z\bigl(\ts\Dep(y)\bigr)=
\id\ox\tau_z\bigl(\ts\De(y)\bigr)\cdot\R(z)
$$
for all $y\in\YN$. For any $r<s$ there is a formal series
$f(z)\in 1+z^{-1}\ts\CC\ts[[z^{-1}]]$ such that $f(z)\ts\R(z)$ acts in
$V_{\la\ru,\,\mu\ru}(h\ru)\ox V_{\la\su,\,\mu\su}(h\su)$ as a rational function
in $z$. This rational function may have zero or pole at $z=0$. Taking
the first non-zero coefficient of the Laurent series of this
function at $z=0$ we get an element
$$
R^{(rs)}\in\End\bigl(\ts V_{\la\ru,\,\mu\ru}\ox V_{\la\su,\,\mu\su}\bigr)
$$
which is an intertwining operator between the {\=$\YN$}-modules
obtained from the tensor product
$V_{\la\ru,\,\mu\ru}(h\ru)\ox V_{\la\su,\,\mu\su}(h\su)$ via
comultiplications $\Dep$ and $\De$.
This intertwining operator is made explicit in our Section 2, it depends
on $h\ru$ and $h\su$ via the difference $h\ru-h\su$.
If this operator is non-invertible for some $r<s$ then
the {\=$\YN$}-module $V$ is reducible.

\proclaim{Conjecture}
The {\=$\YN$}-module $V$ is irreducible if and only if
all the operators $R^{(rs)}$ with $1\le r<s\le n$ are invertible.
\endproclaim

\nt
Our Theorem 3.4 confirms this conjecture when
$\la^{(1)}\lcd\la^{(n)}$ are multiples of any
fundamental weights of $\glN$ while $\mu^{(1)}\lcd\mu^{(n)}$ are empty.
Then Theorem~2.3 describes explicitly the set of differences
$h\ru-h\su\in\ZZ$ where the operator $R^{(rs)}$ is not invertible,
see the remarks after the proof of Theorem 3.4. The proofs of
Theorems 3.3 and 3.4 are based on Proposition 3.1.
It gives suffucient conditions on the parameters $h^{(1)}\lcd h^{(n)}$ 
for cyclicity of the vector $\ze^{(1)}\ox\ldots\ox\ze^{(n)}$
in the {\=$\YN$}-module $V$ with arbitrary $\la^{(1)}\lcd\la^{(n)}$ 
and $\mu^{(1)}\lcd\mu^{(n)}$.

%------------------------------------------------------------------------------

\Section{Elementary modules}
\nt
In this section we consider a class of irreducible
finite-dimensional {\=$\YN$}-modules which arise naturally
from the classical representation theory [C2,O]. Here we also collect
all necessary facts from [MNO\ts,\ts NT] about the Hopf algebra $\YN$.

The \emph{Yangian} of general linear Lie algebra $\glN$
is the associative unital algebra $\YN$ over $\CC$
with the generators $T^{(s)}_{ij}$ where $s=1,2,\ts\ldots$ and
$i,j=1\lcd N$. Defining relations in the algebra $\YN$
can be written for the generating series
$$
T_{ij}(u)=
\de_{ij}+T_{ij}^{(1)}\ts u^{-1}+T_{ij}^{(2)}\ts u^{-2}+\ts\ldots
\nopagebreak
\Tag{1.0101}
$$
in a formal parameter $u$ as follows:
for all indices $i,j,k,l=1\lcd N$ we have
$$
(u-v)\cdot\ts[\ts T_{ij}(u)\ts,\ts T_{kl}(v)\ts]=
T_{kj}(v)\ts T_{il}(u)-T_{kj}(u)\ts T_{il}(v)\,.
\nopagebreak
\Tag{1.1}
$$
Here $v$ is another formal parameter.
\!Let us rewrite these relations \text{in a matrix form.}

Let $E_{ij}\in\EndCN$ be the standard matrix units.
Combine all the series $T_{ij}(u)$ into the single element
$$
T(u)=\sum^N_{i,j=1}\ts E_{ij}\ox T_{ij}(u)
$$
of the algebra $\EndCN\ox\YN\ts[[u^{-1}]]$.
Consider the \emph{Yang {\=$R$}-matrix}
$$
R(u,v)=\id+\sum_{i,j=1}^N\ts\frac{\ts E_{ij}\ox E_{ji}}{u-v}
\in\EndCN^{\ox\ts2}\ts(u,v)\,.
\Tag{1.0}
$$
For any associative unital algebra $\operatorname{X}$ denote
by $\io_s$ its embedding into a finite tensor product
$\operatorname{X}^{\ox\ts n}$
as the $s$th tensor factor:
$$
\io_s(x)=1^{\ox\ts (s-1)}\ox x\ox1^{\ox\ts(n-s)},\quad
x\in\operatorname{X};\qquad s=1\lcd n\,.
$$
Introduce the formal power series with the coefficients
in $\EndCN^{\ox2}\ox\YN$
$$
T_1(u)=\io_1\ox\id\bigl(T(u)\bigr)
\quad\text{and}\quad
T_2(v)=\io_2\ox\id\bigl(T(v)\bigr)\ts.
$$
Then the relations \(1.1) divided by $u-v$
can be rewritten as the single equality
$$
R(u,v)\ox1\cdot T_1(u)\ts T_2(v)=
T_2(v)\ts T_1(u)\cdot R(u,v)\ox1\,.
\Tag{1.2}
$$

Relations \(1.1) imply that for any formal series
$f(u)\in 1+u^{-1}\ts\CC\ts[[u^{-1}]]$ the assignement of generating series
$T_{ij}(u)\mapsto f(u)\cdot T_{ij}(u)$ defines an automorphism of the algebra
$\YN$. We will denote this automorphism by $\of$.

The element $T(u)$ of the algebra $\EndCN\ox\YN\ts[[u^{-1}]]$
is invertible, let us denote
$$
T(u)^{-1}=\widetilde T(u)=\sum^N_{i,j=1}\ts E_{ij}\ox\widetilde{T}_{ij}(u).
$$
Then the relations \(1.2) along with the equality
$
R(u,v)\ts R(-u,-v)=1-(u-v)^{-2}
$
imply that the assignment $T_{ij}(u)\mapsto\widetilde{T}_{ij}(-u)$
determines an automorphism of the algebra $\YN$.
We will denote by this automorphism $\si_N$, it is clearly involutive.

Now consider the elements $E_{ij}$ as generators of the Lie algebra $\glN$.
The algebra $\YN$ contains the enveloping elgebra $\UN$
as a subalgebra: due to \(1.1) the assignment $E_{ji}\mapsto T_{ij}^{(1)}$
defines the embedding. Moreover, there is a homomorphism
$$
\pi_N:\YN\to\UN:\ts T_{ij}(u)\mapsto \de_{ij}+E_{ji}\ts u^{-1}.
\Tag{1.25}
$$
Note that this homomorphism is identical on the subalgebra $\UN$ by definition.
We will fix the Borel subalgebra in $\glN$ generated by the elements $E_{ij}$
with $i\le j$. We will also fix the basis $E_{\ts11}\lcd E_{NN}$ in the
corresponding Cartan subalgebra.

For any non-negative integer $M$ we will fix the standard embedding of
Lie algebras $\glM\to\glMN:\ts E_{ij}\mapsto E_{ij}$.
By \(1.1) there is
an embedding of algebras
$$
\ph:\YN\to\YMN:\ts T_{ij}(u)\mapsto T_{M+i,M+j}(u).
$$
Consider the embedding of the same algebras
$
\ps=\si_{M+N}\ts\ph\,\si_N.
$
The following observation first appeared in [C2]\ts;
for its proof see also [O] and [NT\ts,\ts Section 1].

\proclaim{Proposition 1.1}
Image of the homomorphism $\pi_{M+N}\o\ps:\ts\YN\to\UMN$
commutes with the subalgebra $\UM$ in $\UMN$.
\endproclaim

\nt
This proposition allows us to define a family of {\=$\YN$}-modules which
we will call \emph{elementary}.
For any pair of non-increasing sequences of integers
$$
\la=(\ts\la_1\lcd\la_M,\la_{M+1}\lcd\la_{M+N}\ts)
\quad\text{and}\quad
\mu=(\ts\mu_1\lcd\mu_M\ts)
$$
denote by $V_\lm$ the subspace in the irreducible {\=$\glMN$}-module
of highest weight $\la$ formed by all
singular vectors with respect to $\glM$ of weight $\mu$.
This subspace is preserved by the action of the image of $\pi_{M+N}\o\ps$.
Thus $V_\lm$ becomes a module over the algebra $\YN$.
Relations \(1.1) imply that for any $h\in\CC$
the assignment $T_{ij}(u)\mapsto T_{ij}(u+h)$
determines an automorphism of the algebra $\YN$;
here the series in $(u+h)^{-1}$ should be re-expanded in $u^{-1}$.
We will denote by
$V_\lm(h)$ the {\=$\YN$}-module obtained from $V_\lm$ by pulling back through
this automorphism.

The study of the elementary modules $V_\lm(h)$ has been commenced in [C2]
and continued in [NT]. Let us recall some of these results.
There is a distinguished basis in the
vector space $V_\lm$. It constitutes a part of the \emph{Gelfand-Zetlin basis}
in the irreducible {\=$\glMN$}-module $V_\la$
of the highest weight $\la$, corresponding to the chain of Lie subalgebras
$\gl_1\subset\gl_2\subset\ldots\subset\gl_{M+N}$.
The elements of the latter basis are labelled [GZ]
by the arrays with integral entries
$$
\La=\bigl(\,\la_{mi}\ |\ m=1\lcd M+N\ts;\ i=1\lcd m\,\bigr)
$$
where $\la_{M+N,i}=\la_i$
and $\la_{mi}\ge\la_{m-1,i}\ge\la_{m,i+1}$
for all possible $m$ and $i$. These arrays are called
\emph{Gelfand-Zetlin schemes} of type $\la$. For each scheme $\La$
there is a unique one-dimensional subspace $V_\La\subset V_\la$ which
for every $m=1\lcd M+N-1$ is contained in an irreducible {\=$\gl_m$}-module
of highest weight $(\la_{m1},\la_{m2}\lcd\la_{mm})$.
By choosing a non-zero vector $\xi_\La$ in each subspace $V_\La$ one
obtains a basis in $V_\la$.

The distinguished basis in $V_\lm$ is formed by the vectors $\xi_\La$
where the scheme $\La$ satisfies the condition $\la_{mi}=\mu_i$
for every $m=1\lcd M$. We will denote by $\S_\lm$ the set of all such schemes.
To describe the action of the Yangian $\YN$
on the vectors of this basis in $V_\lm(h)$ explicitly,
it is convenient to use
a set of generators different from $T_{ij}^{(s)}$.
These alternative generators of the algebra $\YN$ are called
the \emph{Drinfeld generators} and can be defined as follows.
Let $\bi=(i_1,\dots,i_k)$ and $\bj=(j_1,\dots,j_k)$ be
any two sequences of integers such that
$$
1\le i_1<\ldots<i_k\le N
\quad\text{and}\quad
1\le j_1<\ldots<j_k\le N.
\Tag{1.999}
$$
Consider the alternating sum over all elements $g$ of
the symmetric group $S_k$
$$
\align
Q^{\vp1}_{\bi\bj}(u)
&\,=\,
\sum_g\
T_{i_1 j_{g(1)}}(u)\ts T_{i_2 j_{g(2)}}(u-1)\ \ldots\ T_{i_k j_{g(k)}}(u-k+1)
\cdot\operatorname{sgn}g
\\
&\,=\,
\sum_g\
T_{i_{g(k)} j_k}(u-k+1)\ \ldots\ T_{i_{g(2)} j_2}(u-1)\ts T_{i_{g(1)} j_1}(u)
\cdot\operatorname{sgn}g
\Tag{1.7777777}
\endalign
$$
where the series
in $(u-1)^{-1},\ts\dots\ts,(u-k+1)^{-1}$ should be
re-expanded in $u^{-1}$. For the proof of the second equality here
see [\ts MNO, Section 2\ts].
For each $k=1\lcd N$ denote
$A_k(u)=Q^{\vp1}_{\bi\bi}(u)$ where $\bi=(1,\ldots,k)$.
Set $A_0(u)=1$.
The series $A_N(u)$ is called the
\emph{quantum determinant} for the Yangian $\YN$.
The next proposition is well known, see [MNO\ts,\,Section 2]
for its detailed proof.

\proclaim{Proposition 1.2}
The coefficients at $u^{-1},\ts u^{-2},\ts\dots$ of the series $A_N(u)$
are free generators for the centre of the algebra $\YN$.
\endproclaim

\nt
This proposition implies that all the coefficients of the series
$A_1(u)\lcd A_N(u)$ pairwise commute.
Further, for each $k=1\lcd N-1$ denote
$$
B_k(u)=Q^{\vp1}_{\bi\bj}(u),\quad
C_k(u)=Q^{\vp1}_{\bj\bi}(u),\quad
D_k(u)=Q^{\vp1}_{\bj\bj}(u)
\Tag{1.99}
$$
where $\bi=(1,\ldots,k)$ and $\bj=(1,\ldots,k-1,k+1)$. The coefficients
of the \text{series} $B_1(u),C_1(u)\lcd B_{N-1}(u),C_{N-1}(u)$ along with
those of
$A_1(u)\lcd A_N(u)$ also generate [D2,\ts Example] the algebra $\YN$.
It is the action of these generators of $\YN$ in $V_\lm(h)$ that
can be calculated explicitly. \text{For any $k=1\lcd N$ \hskip-1pt denote}
$$
\rho_k(u)\,=\,
\prod_{i=1}^{M}\,\ts\frac{u+h+\mu_i-i-k+1}{u+h-i-k+1}
\,\ts\cdot
\prod_{i=1}^{M+k}\ts{\left(u+h-i+1\right)}\ts.
\Tag{1.007}
$$
Regard $\rho_k(u)$ as a formal Laurent series in $u^{-1}$. From now on
we will assume that the space $V_\lm$ is non-zero so that
$S_\lm\neq\varnothing$. Take any scheme $\La$ from $\S_\lm\ts$.

\proclaim{Theorem 1.3}
\!For $k=1\lcd N$ we have equality of formal
\text{Laurent series in $u^{-1}$\hskip-2.1pt}
$$
\rho_k(u)\ts A_k(u)\ts\cdot\,\xi_\La=\xi_\La
\ts\cdot
\prod_{i=1}^{M+k}\left(u+h+\la_{M+k,i}-i+1\right).
\Tag{1.3}
$$
For $k=1\lcd N-1$ formal Laurent series
$\rho_k(u)\ts B_k(u)\cdot\xi_\La$ and $\rho_k(u)\ts C_k(u)\cdot\xi_\La$
in $u^{-1}$ are actually polynomials in $u$
and their degrees are less than $M+k$.
\endproclaim

\nt
This theorem is contained in [NT\ts,\ts Section 2].
Let us now consider the zeroes
$$
\nu_{ki}=i-h-\la_{M+k,i}-1
\,;\qquad
i=1\lcd M+k
\nopagebreak
\Tag{1.4}
$$
of the polynomial in $u$ at the right hand side of \(1.3).
Note that all these $M+k$ zeroes are pairwise distinct since
$\la_{M+k,1}\ge\ldots\ge\la_{M+k,M+k}$ for any Gelfand-Zetlin
scheme $\La$. Therefore the polynomials
$\rho_k(u)\ts B_k(u)\cdot\xi_\La$ and $\rho_k(u)\ts C_k(u)\cdot\xi_\La$
can be determined by their values at the points \(1.4)\ts.
To write down these values one has to make a choice of the vector
$\xi_\La\in V_\La$ for every $\La\in\S_\lm$.
Both these tasks have been performed in [NT\ts,\ts Section 3].
The next twin theorems are weaker than those results of [NT]
but will suffice for our present purposes.
Let the indices $k\in\{1\lcd N-1\}$ and $i\in\{1\lcd M+k\}$ be fixed.
Denote by $\La^\tm$ and $\La^\tp$ the arrays obtained from $\La$ by
decreasing and increasing the {\=$(M+k,i)$}-entry by $1$.

\proclaim{Theorem 1.4}
If $\La^\tm\in\S_\lm\ts$ then
the image $\rho_k(u)\ts B_k(u)\cdot V_\La$ at $u=\nu_{ki}$ is $V_{\La^\tm}$.
If otherwise $\La^\tm\notin\S_\lm\ts$ then this image at $u=\nu_{ki}$ is zero.
\endproclaim

\proclaim{Theorem 1.5}
If $\La^\tp\in\S_\lm\ts$ then
the image $\rho_k(u)\ts C_k(u)\cdot V_\La$ at $u=\nu_{ki}$ is $V_{\La^\tp}$.
If otherwise $\La^\tp\notin\S_\lm\ts$ then this image at $u=\nu_{ki}$ is zero.
\endproclaim

\nt
The above three theorems imply that the elementary {\=$\YN$}-module $V_\lm(h)$
is irreducible. Let us now point out the place of elementary modules
in the family of all irreducible finite-dimensional {\=$\YN$}-modules.
Let $V$ be any module from the latter family.
A non-zero vector $\ze$ in any {\=$\YN$}-module is called \emph{singular}
if it is annihilated by all the coefficients of the series
$C_1(u)\lcd C_{N-1}(u)$. The vector $\ze\in V$ is then unique
up to a scalar multiplier and is an eigenvector for the coefficients of the
series $A_1(u)\lcd A_N(u)$; see [D2,\ts Theorem 2]. Moreover, then
$$
\frac
{A_{k+1}(u)\ts A_{k-1}(u-1)}
{A_k(u)\ts A_k(u-1)}\cdot\ze
=
\frac{P_k(u-1)}{P_k(u)}\cdot\ze\,;
\qquad
k=1\lcd N-1
\Tag{1.9999}
$$
for certain monic polynomials $P_1(u)\lcd P_{N-1}(u)$ with coefficients
in $\CC$.
These $N-1$ polynomials are called the \emph{Drinfeld polynomials} of the
module $V$.
Every collection of $N-1$ monic polynomials arises in this way. The modules
with the
same Drinfeld polynomials may differ only by an automorphism of the algebra
$\YN$ of the form $\of$.
Now consider the scheme $\La^{\hskip-1pt\circ}\in\S_\lm$ with the entries
$$
\la_{mi}^\circ=
\cases
\mu_i
&\quad\text{if}\quad m\le M,
\\
\operatorname{min}(\ts\la_i,\mu_{i-m+M}\ts)
&\quad\text{if}\quad m>M\ \ \text{and}\ \ i>m-M,
\\
\la_i
&\quad\text{if}\quad m>M\ \ \text{and}\ \ i\le m-M.
\endcases
$$
Observe that $\la_{mi}^\circ\ge\la_{mi}$ for every scheme $\La\in\S_\lm$.
Therefore by Theorem 1.5 the vector $\xi_{\La^{\hskip-1pt\circ}}\in V_\lm(h)$
is singular. Theorem 1.3 then allows to find the Drinfeld polynomials of the
module $V_\lm(h)$. We again refer to [\ts NT\ts,\, Section 2\ts] for details.

Here we will only formulate the answer.
The assumption $S_\lm\neq\varnothing$ implies that
$$
\la_1\ge\mu_1
\ts\lcd
\la_M\ge\mu_M
\quad\text{and}\quad
\mu_M\ge\la_{M+N}\,.
\Tag{1.5}
$$
Consider the \emph{skew Young diagram} $\lsm$. This is the set of pairs
$$
\{\,(i,j)\in\ZZ^2\ |\ 1\le i\le M+N,\ \la_i\ge j>\mu_i\,\}
$$
where for any $i>M$ we write $\mu_i=\la_{M+N}$.
Employ usual graphic representation of a diagram:
the point $(i,j)\in\ZZ^2$
is represented by the unit box on the plane $\Bbb R^2$ with the centre
$(i,j)$;
the coordinates $i$ and $j$ on $\Bbb R^2$ increasing from top to bottom
and from left to
right respectively. The \emph{content} of the box corresponding to $(i,j)$
is the
difference $c=j-i$. Here is the diagram corresponding to
$\la=(\ts5,\ts5,\ts3,\ts2,\ts0,\mi2)$ and $\mu=(\ts3,\ts2,\ts2,\ts1)$;
we indicate the content of the bottom box for every column.
\medskip
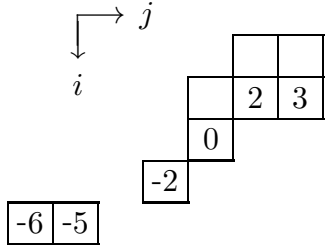
\vbox{
$$
\kern23pt\longrightarrow\,j\mv
$$
\vglue-26.2pt
$$
\vert\mv
$$
\vglue-28pt
$$
\bigr\downarrow\mv\kern-.1pt
$$
\vglue-16pt
$$
i\mv
$$
\vglue-40pt
$$
\phantom{\bx}
\phantom{\bx}
\phantom{\bx}
\phantom{\bx}
\phantom{\bx}
{\bx}
{\bx}
$$
\vglue-17.8pt
$$
\phantom{\bx}
\phantom{\bx}
\phantom{\bx}
\phantom{\bx}
{\bx}
{\bx}
{\bx}
$$
\vglue-17.8pt
$$
\phantom{\bx}
\phantom{\bx}
\phantom{\bx}
\phantom{\bx}
{\bx}
\phantom{\bx}
\phantom{\bx}
$$
\vglue-17.8pt
$$
\phantom{\bx}
\phantom{\bx}
\phantom{\bx}
{\bx}
\phantom{\bx}
\phantom{\bx}
\phantom{\bx}
$$
\vglue-17.8pt
$$
{\bx}
{\bx}
\phantom{\bx}
\phantom{\bx}
\phantom{\bx}
\phantom{\bx}
\phantom{\bx}
$$
\vglue-70.8pt
$$
\kern71pt
2
\kern9pt
3
$$
\vglue-18pt
$$
\kern29pt
0
$$
\vglue-17.8pt
$$
\mi2
$$
\vglue-17.8pt
$$
\phantom{\mi7}
\kern6pt
\mi6
\kern6pt
\mi5
\kern85pt
$$
}
\medskip
\nt
The condition $\S_\lm\neq\varnothing$ is then equivalent to \(1.5) along with
the requirement that any column of the skew diagram $\lsm$ has at most
$N$ boxes; see \text{[\ts M\ts,\ts Section I.5\ts].}

\proclaim{Proposition 1.6}
The Drinfeld polynomials of the {\=$\YN$}-module $V_\lm(h)$ are
$$
P_k(u)\ts=\ts\prod_c\,\,(u+h+c)\,;
\qquad
k=1\lcd N-1
\nopagebreak
$$
where the product is taken over the contents of the bottom boxes in
the columns of height $k$ in the skew Young diagram $\lsm$.
\endproclaim

\nt
Let us equip the algebra $\YN$ with the {\=$\ZZ$}-grading $\ts\deg\!$
determined by
$$
\deg T_{ij}^{(s)}=i-j;\qquad s=0,1,2,\dots\ts;
\Tag{1.9}
$$
see defining relations \(1.1).
We will extend this grading to the algebra
$\YN[[u^{-1}]]$ by assuming that $\deg u^{-1}=0$.
Then $\deg A_k(u)=0$ while by the definition \(1.99)
$$
\deg B_k(u)=-1\ts,
\quad
\deg C_k(u)=1\ts.
$$

The vector space $V_\lm$ has a natural {\=$\ZZ$}-grading
by eigenvalues of the action of
$$
(M+N)\ts E_{11}+(M+N-1)\ts E_{22}+\ldots+E_{M+N,M+N}\,\in\,\glMN\ts.
$$
Any subspace $V_\La$ in $V_\lm$
is homogeneous and its degree equals the sum of $\la_{mi}$ over all
indices $m=1\lcd M+N$ and $i=1\lcd m$.
The above three theorems show that the action of the algebra $\YN$
in the module $V_\lm(h)$ is graded. The singular vector
$\xi_{\La^{\hskip-1pt\circ}}\in V_\lm(h)$ then has the maximal degree.
Hence it is an eigenvector for all the coefficients of the series
$T_{11}(u)\lcd T_{NN}(u)$.
Let $R_1(u)\lcd R_N(u)$ be the corresponding eigenvalues.
Let $\operatorname{J}$ be the left ideal in $\YN\ts[[u^{-1}]]$
generated by the elements of positive {\=$\ZZ$}-degrees.
By definition $A_k(u)$ equals
$$
T_{11}(u)\,T_{22}(u-1)\ts\ldots\ts T_{kk}(u-k+1)
$$
plus certain elements from the ideal $\operatorname{J}$. Therefore
from \(1.9999) we get the equalities
$$
\frac{P_k(u-1)}{P_k(u)}=
\frac{R_{k+1}(u-k)}{R_k(u-k)}\,;
\qquad
k=1\lcd N-1\,.
\Tag{1.999999999}
$$

Due to the defining relations \(1.1) the assignment
$T_{ij}(u)\mapsto T_{ji}(u)$ determines an anti-automorphism of the
algebra $\YN$. Denote by $\theta$ this anti-automorphism, it
is obviously involutive. Now for any finite-dimensional {\=$\YN$}-module $W$
define its \emph{dual} module $W^\ast$ as the vector space dual to $W$ where
$
\langle\,y\cdot\xi^\ast,\xi\,\rangle=
\langle\,\xi^\ast,\theta(y)\cdot\xi\,\rangle
$
for any $\xi\in W$, $\xi^\ast\in W^\ast$ and $y\in\YN$.
We will need the following simple fact.

\proclaim{Proposition 1.7}
Any elementary {\=$\YN$}-module $V_\lm(h)$ is self-dual.
\endproclaim

\demo{Proof}
Let $\ze\in V_\lm(h)$ be a singular vector. Its spans the subspace of
the maximal {\=$\ZZ$}-degree, so there is an element $\zes\in V_\lm(h)^\ast$
with $\langle\,\zes,\ze\,\rangle=1$ and $\langle\,\zes,\xi\,\rangle=0$
for any vector $\xi\in V_\lm(h)$ with non-maximal degree.
But thanks to \(1.7777777) we have
$$
\theta\bigl(A_k(u)\bigr)=A_k(u)\,,
\quad
\theta\bigl(B_k(u)\bigr)=C_k(u)\,,
\quad
\theta\bigl(C_k(u)\bigr)=B_k(u)
$$
for all possible indices $k$.
Therefore the vector $\zes$ is singular in $V_\lm(h)^\ast$
and the eigenvalues of $A_1(u)\lcd A_N(u)$ on this vector are the same as
on the vector $\ze$ respectively.
So the {\=$\YN$}-modules $V_\lm(h)^\ast$ and $V_\lm(h)$ are equivalent.
\enddemos

\nt
There is a natural Hopf algebra structure on
the $\YN$. The antipode is defined by
the assignment of generating series
$T_{ij}(u)\mapsto\tilde T_{ij}(u)$ while the comultiplication
$\YN\to\YN^{\ox2}$ is defined by the assignement
$$
T_{ij}(u)\,\mapsto\,\sum_{k=1}^N\,T_{ik}(u)\ox T_{kj}(u).
\Tag{1.6}
$$
Here and in what follows we take tensor products
of the elements of the algebra $\YN[[u^{-1}]]$
over its subalgebra $\CC\ts[[u^{-1}]]$.
We will also use the comultiplication $\Dep$ on the algebra $\YN$ obtained
by composing $\De$ with the transposition of tensor factors in $\YN^{\ox2}$.
Observe that by the definition \(1.6) we have
$$
\De\o\theta=(\ts\theta\ox\theta\ts)\o\Dep\,.
\Tag{1.666666}
$$
We will consider the images of the Drinfeld generators for
the algebra $\YN$ with respect to the {\=$n$}-fold comultiplication
$$
\Den:\YN\to\YN^{\ox n}.
\Tag{1.7}
$$
For this purpose we we will employ the following easy
result from [NT, Section 1].
Let $\bi$ and $\bj$ be any two sequences of indices
satisfying the condition \(1.999).

\proclaim{Proposition 1.8}
We have the equality
$$
\Den\bigl(\ts Q_{\bi\bj}(u)\bigr)\ =
\sum_{\bk^{(1)},\bk^{(2)}\lcd\bk^{(n-1)}}\ts
Q_{\bi\bk^{(1)}}(u)\ox
Q_{\bk^{(1)}\bk^{(2)}}(u)\ox\ldots\ox
Q_{\bk^{(n-1)}\bj}(u)
$$
\line{where $\bk^{(1)}\!,\bk^{(2)}\!\lcd\bk^{(n-1)}\!\!$ are
increasing sequences of integers $1\lcd N\!$ of length $\!k.$\!\!\!}
\endproclaim

\nt
With the {\=$\ZZ$}-grading $\deg\!$ on $\YN$, the algebras $\YN^{\ox n}$
and $\YN^{\ox n}[[u^{-1}]]$ acquire grading by the group $\ZZ^n$.
In the next section we will give an alternative realization of the
elementary {\=$\YN$}-module $V_\lm(h)$ using comultiplication \(1.7).

%------------------------------------------------------------------------------

\bigskip
\Section{Intertwining operators}
\nt
In this section we study intertwining operators between the tensor
products of two {\=$\YN$}-modules of the form $V_\lm(h)$ defined via the
comultiplications $\De$ and $\Dep$ on the algebra $\YN$.
We use the explicit realization [C2] of the module $V_\lm(h)$.

Consider first the \emph{vector} {\=$\YN$}-module $V(h)$.
This is the elementary module $V_\lm(h)$ with $M=0$ and $\la=(1,0\lcd 0)$.
So the space of this module is $\CC^N$ and the action of the algebra $\YN$ is
defined by $T_{ij}^{(s)}\mapsto E_{ji}\cdot h^{\ts s-1}$ for any $s\ge1$.
Note
that due to \(1.0) this action can be then determined by the single assignment
$$
\EndCN\ox\YN\ts[[u^{-1}]]\to\EndCN^{\ox2}\ts[[u^{-1}]]:\
T(u)\mapsto R(u,-h)\,.
\Tag{2.0}
$$

Now take any non-zero elementary {\=$\YN$}-module $V_\lm(h)$.
Let $n$ be the total number of boxes in the skew Young diagram $\la/\mu$.
Consider the \emph{row tableau} of shape $\lsm$ obtained by filling
the boxes of $\la$ with $1\lcd n$ in the natural order,
that is by rows downwards, from the left to right in every row.
\text{Denote by $\Om$ this tableau.}
Here is the row tableau corresponding to
$\la=(\ts5,\ts5,\ts3,\ts2,\ts0,\mi2)$ and $\mu=(\ts3,\ts2,\ts2,\ts1)$:
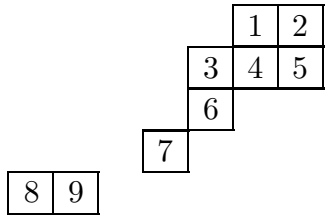
\vsk>
\vbox{
$$
\phantom{\bx}
\phantom{\bx}
\phantom{\bx}
\phantom{\bx}
\phantom{\bx}
{\bx}
{\bx}
$$
\vglue-17.8pt
$$
\phantom{\bx}
\phantom{\bx}
\phantom{\bx}
\phantom{\bx}
{\bx}
{\bx}
{\bx}
$$
\vglue-17.8pt
$$
\phantom{\bx}
\phantom{\bx}
\phantom{\bx}
\phantom{\bx}
{\bx}
\phantom{\bx}
\phantom{\bx}
$$
\vglue-17.8pt
$$
\phantom{\bx}
\phantom{\bx}
\phantom{\bx}
{\bx}
\phantom{\bx}
\phantom{\bx}
\phantom{\bx}
$$
\vglue-17.8pt
$$
{\bx}
{\bx}
\phantom{\bx}
\phantom{\bx}
\phantom{\bx}
\phantom{\bx}
\phantom{\bx}
$$
\vglue-84pt
$$
\kern71pt
1
\kern9pt
2
$$
\vglue-17.8pt
$$
\kern57pt
3
\kern9pt
4
\kern9pt
5
$$
\vglue-18pt
$$
\kern29pt
6
$$
\vglue-17.8pt
$$
7
$$
\vglue-17.8pt
$$
8
\kern9pt
9
\kern70pt
$$
}
\medskip
\nt
Denote by $S_\lsm$ and $T_\lsm$ the subgroups in the symmetric group
$S_n$ preserving the sets of numbers appearing respectively in every row
and in every column of $\Om$.

The symmetric group $S_n$ acts in the space
$(\CC^N)^{\ox n}$ by permutations of the tensor factors.
Let $P_\lsm$ and $Q_\lsm$ be the elements of $\EndCN^{\ox n}$
corresponding to the sums
$$
\sum_{g\in S_\lsm}g
\ \ \quad\text{and}\quad
\sum_{g\in T_\lsm}g\cdot\operatorname{sgn}g
$$
in $\CC\cdot S_n$.
The product $Y_\lsm=P_\lsm\ts Q_\lsm$
is the \emph{Young symmetrizer}
corresponding to the tableau $\Om$.
Let $c_1\lcd c_n$ be the contents of the boxes of $\la/\mu$
occupied respectively by the numbers $1\lcd n$ in the row tableau $\Om$.
Then consider the {\=$\YN$}-module obtained from the tensor product
$V(c_1+h)\ox\ldots\ox V(c_n+h)$ via the comultiplication \(1.7).

\proclaim{Proposition 2.1}
Action of the algebra $\YN$ in \text{$V(c_1+h)\ox\ldots\ox V(c_n+h)$}
preserves the image of the operator $Y_{\lsm}$.
The module $V_\lm(h)$ can be obtained from this image
by pulling back through an automorhism of the form $\of$.
\endproclaim

\demo{Proof}
Consider the {\=$\YN$}-module $V(c_1+h)\ox\ldots\ox V(c_n+h)$.
By \(1.6) and \(2.0) the action of $\YN$ in this module
can be determined by the assignment
$$
\gather
\EndCN\ox\YN\ts[[u^{-1}]]\,\to\,\EndCN^{\ox(n+1)}\ts[[u^{-1}]]:\
\\
T(u)\,\mapsto\,R_{12}(u,-c_1\!-\!h)\ts\ldots\ts R_{1,n+1}(u,-c_n\!-\!h)\,.
\Tag{2.1}
\endgather
$$
Here we use the standard notation:
for any $1\le p<q\le m$ and $X\in\EndCN^{\ox2}$
we write
$$
X_{pq}=\io_p\ox\io_q\,(X)\in\EndCN^{\ox m}\ts.
$$

Introduce the \emph{inverse column tableau} of shape $\lsm$.
This tableau is obtained by filling
the boxes of $\lsm$ with $1\lcd n$ by columns from the right to the left,
upwards in every column. Denote it by $\Om^\ast$.
Here is the tableau $\Om^\ast$ corresponding to the same sequences
$\la=(\ts5,\ts5,\ts3,\ts2,\ts0,\mi2)$ and $\mu=(\ts3,\ts2,\ts2,\ts1)$
as before:
\vsk>
\vbox{
$$
\phantom{\bx}
\phantom{\bx}
\phantom{\bx}
\phantom{\bx}
\phantom{\bx}
{\bx}
{\bx}
$$
\vglue-17.8pt
$$
\phantom{\bx}
\phantom{\bx}
\phantom{\bx}
\phantom{\bx}
{\bx}
{\bx}
{\bx}
$$
\vglue-17.8pt
$$
\phantom{\bx}
\phantom{\bx}
\phantom{\bx}
\phantom{\bx}
{\bx}
\phantom{\bx}
\phantom{\bx}
$$
\vglue-17.8pt
$$
\phantom{\bx}
\phantom{\bx}
\phantom{\bx}
{\bx}
\phantom{\bx}
\phantom{\bx}
\phantom{\bx}
$$
\vglue-17.8pt
$$
{\bx}
{\bx}
\phantom{\bx}
\phantom{\bx}
\phantom{\bx}
\phantom{\bx}
\phantom{\bx}
$$
\vglue-84pt
$$
\kern71pt
4
\kern9pt
2
$$
\vglue-17.8pt
$$
\kern57pt
6
\kern9pt
3
\kern9pt
1
$$
\vglue-18pt
$$
\kern29pt
5
$$
\vglue-17.8pt
$$
7
$$
\vglue-17.8pt
$$
9
\kern9pt
8
\kern70pt
$$
}
\medskip
\nt
Let the symmetric group $S_n$ act on the entries of the tableau $\Om^\ast$.
Consider the permutation $g\in S_n$ such that $g:\Om^\ast\mapsto\Om$.
Then [\ts C1,\ts Theorem\,1\ts] provides the equality in the algebra
$\EndCN^{\ox(n+1)}\ts[[u^{-1}]]$
$$
\gather
R_{12}(u,-c_1\!-\!h)
\ts\ldots\ts
R_{1,n+1}(u,-c_n\!-\!h)
\cdot(1\ox Y_\lsm)=
\Tag{2.2}
\\
(1\ox Y_\lsm)\cdot
R_{1,g(1)+1}(u,-c_{g(1)}\!\!-\!h)
\ts\ldots\ts
R_{1,g(n)+1}(u,-c_{g(n)}\!\!-\!h)\,.
\endgather
$$
Along with \(2.1) this equality yeilds the first statement of
Proposition 2.1. We will denote by $V_\lsm(h)$ the image
of the operator $Y_\lsm$ regarded as {\=$\YN$}-module.

Further, the symmetric group $S_n$ acts on the algebra $\YN^{\ox n}$
by permutations
of the tensor factors. Consider the {\=$\YN$}-module $U$ obtained from
the product $V(c_1+h)\ox\ldots\ox V(c_n+h)$ by composing
the comultiplication $\Den$ with the above permutation $g$.
The action of the algebra $\YN$ in $U$ can be then described by
$$
\gather
\EndCN\ox\YN\ts[[u^{-1}]]\,\to\,\EndCN^{\ox(n+1)}\ts[[u^{-1}]]:\
\\
T(u)\,\mapsto\,
R_{1,g(1)+1}(u,-c_{g(1)}\!\!-\!h)
\ts\ldots\ts
R_{1,g(n)+1}(u,-c_{g(n)}\!\!-\!h)\,.
\endgather
$$
Equality \(2.2) also shows that this action
preserves the kernel of the operator $Y_\lsm$.
Moreover, by \(2.2) the {\=$\YN$}-module $V_{\la/\mu}(h)$ is equivalent
to the quotient of $U$ by this kernel. We will show that this quotient is
irreducible and has the same Drinfeld polynomials as
the elementary {\=$\YN$}-module $V_\lm(h)$.

First consider the case when the skew Young diagram $\la/\mu$ consists of
one column only. By our assumption the length $n$ of this column does not
exceed $N$. Then we have $c_p=c_1-p+1$ and $g(p)=n-p+1$ for each $p=1\lcd n$
while $Y_\lsm$ is the operator of antisymmetrization
in $(\CC^N)^{\ox n}$. Further, then for
$$
P=\sum^N_{i,j=1}\ts E_{ij}\ox E_{ji}\,\in\,\EndCN^{\ox2}
$$
we have the equality of rational functions in $u$
valued in the algebra $\EndCN^{\ox(n+1)}\!$
$$
\gather
(1\ox Y_\lsm)\cdot
R_{1,n+1}(u,-c_n\!\!-\!h)
\ts\ldots\ts
R_{12}  (u,-c_1\!\!-\!h)
\\
=(1\ox Y_\lsm)\cdot
\Bigl(1+\frac{P_{12}+\ldots+P_{1,n+1}}{u+h+c_1}\Bigr)\,;
\endgather
\nopagebreak
$$
see for instance [\ts N,\,Proposition 2.12\ts].
This equality along with \(1.25)
shows that the quotient of $U$ by $\ker Y_\lsm$
is equivalent to the {\=$\YN$}-module
$V_{\nu,\raise1pt\hbox{$\hskip-1pt{}_\varnothing$}}(c_1+h)$ where $\nu$ is
the $n$th fundamental {\=$\glN$}-weight $(1\lcd 1,0\lcd 0)$.
Let $e_1\lcd e_N$ be the standard basis in $\CC^N$.
Due to Proposition 1.6 every Drinfeld polynomial of the {\=$\YN$}-module
$V_{\nu,\raise1pt\hbox{$\hskip-1pt{}_\varnothing$}}(c_1+h)$ is $1$,
except for the $n$th when $n<N$. \text{Then the} $n$th polynomial is
$$
u+c_1+h+1-n=u+h+c_n\,.
$$
The Drinfeld polynomials of the elementary module $V_\lm(h)$ are the same.
So we obtain the required statement in the one-column case.
Observe that in this case the image of the vector $e_1\ox\ldots\ox e_n\in U$
in the quotient by $\ker Y_\lsm$ is singular.

Now consider the case of arbitrary $\la$ and $\mu$. Let us denote by $\S_\lsm$
the set of \emph{semi-standard tableaux} of shape $\lsm$ with entries
$1\lcd N$. These tableaux are the functions $\kappa:\{1\lcd n\}\to\{1\lcd N\}$
such that $\ka(p)<\ka(q)$ or $\ka(p)\le\ka(q)$ if the numbers $p<q$
appear respectively in the same column or \text{the same row of $\Om$.}
Every tableau $\ka\in\S_\lsm$ determines a vector
$e_{\ka(1)}\ox\ldots\ox e_{\ka(n)}\in(\CC^N)^{\ox n}$.
The images of all these vectors in the quotient of $(\CC^N)^{\ox n}$
by $\ker Y_\lsm$ form a basis in this quotient; see for instance
[\ts JK,\ts Section 7.2\ts].
Moreover, the sets $\S_\lsm$ and $\S_\lm$ are of the same cardinality.
Therefore it suffices to point out a singular vector $\ze$ in the quotient
by $\ker Y_\lsm$ of the {\=$\YN$}-module $U$ such that the equaities \(1.9999)
hold for the Drinfeld polynomials $P_1(u)\lcd P_{N-1}(u)$ of
the module $V_\lm(h)$.
In particular, we will then obtain that this quotient is irreducible.

Let $m$ be the number of columns in the skew diagram $\lsm$.
By our assumption the length of any column does not exceed $N$.
Let $\ka^\circ(p)\in\{1\lcd N\}$ be the depth of the box of the tableau
$\Om$ with the number $p$ in its column. Then $\ka^\circ\in\S_\lsm$.
Arguments already used in the one-column case show that the action
of the algebra $\YN$ in $U$ preserves the kernel of the operator $Q_\lsm$.
The image $\eta$
of the vector $e_{\ka^\circ(1)}\ox\ldots\ox e_{\ka^\circ(n)}\in U$
in the quotient by $\ker Q_\lsm$ is singular.
This follows from Proposition 1.8 applied to the number $m$ instead of $n$,
see also \(1.9). Moreover, the same proposition shows that for any $k=1\lcd N$
the coproduct $\De^{\hskip-1pt(m)}(A_k(u))$ equals $A_k(u)^{\ox m}$ plus
terms with degrees in $\ZZ^m$ containing at least one positive component.
Using the results of the one-column case, we then obtain the equalities
$$
\frac
{A_{k+1}(u)\ts A_{k-1}(u-1)}
{A_k(u)\ts A_k(u-1)}\cdot\eta
\,=\,\eta\,\cdot\,\prod_c\,\frac{u+h+c-1}{u+h+c}\,\,;
\quad\ k=1\lcd N-1
\Tag{2.9999}
$$
where the product is taken over the contents of the bottom boxes in
the columns of height $k$ in the skew Young diagram $\lsm$.
Now let $\ze$ be the image of the vector
$e_{\ka^\circ(1)}\ox\ldots\ox e_{\ka^\circ(n)}\in U$
in the quotient by $\ker Y_\lsm$. Since
$\ker Q_\lsm\subset\ker Y_\lsm$, the vector $\ze$ is singular in this
quotient {\=$\YN$}-module.
Moreover, the equalities \(2.9999) along with Proposition 1.6 imply \(1.9999)
for this quotient and for the Drinfeld polynomials $P_1(u)\lcd P_{N-1}(u)$ of
the elementary module $V_\lm(h)$.
\enddemos

\nt
Now fix any two skew Young diagrams $\al$ and $\be$.
Let $m$ and $n$ be the numbers of boxes in $\al$ and $\be$ respectively.
Let $z$ be a complex parameter \text{as well as $h$.}
Consider the irreducible {\=$\YN$}-modules $V_\al(h)$ and $V_\be(z)$.
Let $a_1\lcd a_m$ and $b_1\lcd b_n$ be the contents of the boxes
of $\al$ and $\be$ occupied by the numbers $1\lcd m$ and $1\lcd n$ in the
corresponding row tableaux. Then $V_\al(h)$ and $V_\be(z)$ are the submodules
in $V(a_1+h)\ox\ldots\ox V(a_m+h)$ and $V(b_1+z)\ox\ldots\ox V(b_n+z)$
defined as the images of Young symmetrizers $Y_\al$ and $Y_\be$
in $(\CC^N)^{\ox m}$ and $(\CC^N)^{\ox n}$
respectively; see the proof of Proposition 2.1.

We assume that the {\=$\YN$}-modules $V_\al(h)$ and $V_\be(z)$
are both non-zero. So the length of any column in $\al$ and $\be$
does not exceed $N$.
Note that if $m=n$ while $a_k=b_k+c$ for each $k=1\lcd m$
and the same integer $c$ then $V_\al(h)=V_\be(h+c)$.

Equip the set all of pairs $(k,l)$ where $k=1\lcd m$ and $l=1\lcd n$ with
the following ordering: the pair $(i,j)$ precedes $(k,l)$ if $i>k$,
or if $i=k$ but $j<l$.
Using this ordering introduce the rational function in $\!h,z\!$
valued in $\EndCN^{\ox(m+n)}\!$
$$
\prod_{(k,l)}^\rightarrow\ R_{k,m+l}\,(-a_k\!-\!h,-b_l\!-\!z)
\,\cdot\,Y_\al\ox Y_\be
\Tag{2.3}
$$
where $Y_\al$ acts non-trivially only on the first $m$ tensor factors in
$(\CC^N)^{\ox(m+n)}$ while $Y_\be$ acts only on the last $n$ tensor factors.
By the definition \(1.0) this function depends only on the difference $h-z$.
Due to \(2.2) the expression \(2.3) is divisible by $1\ox Y_\be$
also on the left. Furthermore, since $R(u,v)\ts R(v,u)=1-(u-v)^{-2}$
the equality \(2.2) in the algebra $\EndCN^{\ox(n+1)}\ts[[u^{-1}]]$
can be rewritten as
$$
\gather
R_{n,n+1}(-c_n\!-\!h\ts,u)
\ts\ldots\ts
R_{1,n+1}(-c_1\!-\!h\ts,u)
\cdot(Y_\lsm\ox1)=
\\
(Y_\lsm\ox1)\cdot
R_{g(n),n+1}(-c_{g(n)}\!\!-\!h\ts,u)
\ts\ldots\ts
R_{g(1),n+1}(-c_{g(1)}\!\!-\!h\ts,u)\,.
\endgather
$$
Hence the expression \(2.3) is divisible by $Y_\al\ox1$
on the left. Thus (2.5) determines a rational function in $h-z$ valued
in $\End\bigl(\ts\im Y_\al\ox\im Y_\be\bigr)$.
Denote by $R_{\al\be}(h)$ the first non-zero coefficient
of the Laurent series in $z$ of this function at $z=0$.

Now let $W$ and $W^\prime$ be the {\=$\YN$}-modules obtained from the tensor
product $V_\al(h)\ox V_\be(0)$ via the comultiplications $\De$ and $\Dep$
respectively. Then the element
$R_{\al\be}(h)\in\End\bigl(\ts\im Y_\al\ox\im Y_\be\bigr)$
admits the following interpretation.

\proclaim{Proposition 2.2}
The coefficient $R_{\al\be}(h)$
is an intertwining operator $W^\prime\to W$.
\endproclaim

\demo{Proof}
Let us
denote by $U$ and $U^\prime$ the {\=$\YN$}-modules obtained respectively
via the comultiplications $\De$ and $\Dep$ from the tensor product of
the {\=$\YN$}-modules $V(a_1+h)\ox\ldots\ox V(a_m+h)$ and
$V(b_1+z)\ox\ldots\ox V(b_n+z)$. The action of the algebra $\YN$
in the module $U$ can be determined by the assignment
$$
\gather
\EndCN\ox\YN\ts[[u^{-1}]]\,\to\,\EndCN^{\ox(m+n+1)}\ts[[u^{-1}]]:\
\\
T(u)\,\mapsto\,
R_{12}(u,-a_1\!-\!h)\ts\ldots\ts R_{1,m+1}(u,-a_m\!-\!h)\,\x
\\
\qquad\qquad
R_{1,m+2}(u,-b_1\!-\!z)\ts\ldots\ts R_{1,m+n+1}(u,-b_n\!-\!z)
\,;
\endgather
$$
cf. \(2.1). The action of $\YN$ in $U^\prime$
can be then determined by the assignment
$$
\gather
\EndCN\ox\YN\ts[[u^{-1}]]\,\to\,\EndCN^{\ox(m+n+1)}\ts[[u^{-1}]]:\
\\
T(u)\,\mapsto\,
R_{1,m+2}(u,-b_1\!-\!z)\ts\ldots\ts R_{1,m+n+1}(u,-b_n\!-\!z)\,\x
\\
\qquad\qquad
R_{12}(u,-a_1\!-\!h)\ts\ldots\ts R_{1,m+1}(u,-a_m\!-\!h)\,.
\endgather
$$
By applying repeatedly the \emph{Yang-Baxter equation} in $\EndCN^{\ox3}(u,v,w)$
$$
R_{12}(u,v)\,R_{13}(u,w)\,R_{23}(v,w)=
R_{23}(v,w)\,R_{13}(u,w)\,R_{12}(u,v)
$$
we obtain the equality of rational functions in $h\ts,z$ valued
in $\EndCN^{\ox(m+n+1)}$
$$
\gather
\prod_k^\rightarrow\,R_{1,k+1}(u,-a_k\!-\!h)\,\cdot\,
\prod_l^\rightarrow\,R_{1,m+l+1}(u,-b_l\!-\!z)\,\,\x
\qquad\qquad
\\
\ \qquad
\prod_{(i,j)}^\rightarrow\ R_{k+1,m+l+1}\,(-a_k\!-\!h,-b_l\!-\!z)\,\cdot\,
1\ox Y_\al\ox Y_\be\,\,=
\\
\qquad\qquad
\prod_{(k,l)}^\rightarrow\ R_{k+1,m+l+1}\,(-a_k\!-\!h,-b_l\!-\!z)\,\,\x
\\
\prod_l^\rightarrow\,R_{1,m+l+1}(u,-b_l\!-\!z)\,\cdot\,
\prod_k^\rightarrow\,R_{1,k+1}(u,-a_k\!-\!h)\,\cdot\,
1\ox Y_\al\ox Y_\be\,.
\endgather
$$
Here the index $k$ runs through the set
$\!\{1\lcd m\}\!$ while $l$ runs through $\!\{1\lcd n\}\!$.
Note that the expression in the last line of this equality is
divisible by $1\ox Y_\al\ox Y_\be$ also on the left.
Let $d$ be the degree of the first non-zero term in the Laurent series
of the rational function \(2.3) in $z$ at $z=0$.
Dividing the above equality by $z^d$ and then tending $z\to0$,
we obtain Proposition 2.2.
\enddemos

\nt
By the definition \(2.3) the operator $R_{\al\be}(h)$ in
$\im Y_\al\ox\ts\im Y_\be$ is invertible for every $h\in\CC\setminus\ZZ$.
Now consider the following special situation. Suppose that
the diagram $\be$ is a usual Young diagram. Thus
for certain integers $\be_1\ge\ldots\ge\be_N\ge0$ we have
$$
\be=\{\,(i,j)\in\ZZ^2\ |\ 1\le i\le N,\
0<j\le\be_i\,\}\,.
$$
Note that in this case $b_1\!=0$.
Further, suppose that $\al$ is a \emph{reversed} Young diagram:
$$
\al=\{\,(i,j)\in\ZZ^2\ |\ 1\le i\le N,\
N-\al_{N-i+1}< j\le N\,\}
\nopagebreak
$$
for certain integers $\al_1\ge\ldots\ge\al_N\ge0$.
Note that then $a_m=0$. For example,
here for $N=3$ we show the reversed and the
usual Young diagrams corresponding to
$(\al_1,\al_2,\al_3)=(\be_1,\be_2,\be_3)=(4,3,1)$.
We have indicated the contents $a_8=b_1=0$.
\smallskip
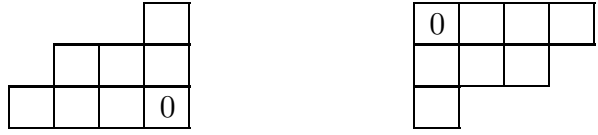
\vsk>
\vbox{
$$
\phantom{\bx}
\phantom{\bx}
\phantom{\bx}
{\bx}
\phantom{\bx}
\phantom{\bx}
\phantom{\bx}
\phantom{\bx}
\phantom{\bx}
{\bx}
{\bx}
{\bx}
{\bx}
$$
\vglue-17.8pt
$$
\phantom{\bx}
{\bx}
{\bx}
{\bx}
\phantom{\bx}
\phantom{\bx}
\phantom{\bx}
\phantom{\bx}
\phantom{\bx}
{\bx}
{\bx}
{\bx}
\phantom{\bx}
$$
\vglue-17.8pt
$$
{\bx}
{\bx}
{\bx}
{\bx}
\phantom{\bx}
\phantom{\bx}
\phantom{\bx}
\phantom{\bx}
\phantom{\bx}
{\bx}
\phantom{\bx}
\phantom{\bx}
\phantom{\bx}
$$
\vglue-57.5pt
$$
\kern86pt{0}
$$
\vglue-4.8pt
$$
{0}\kern84pt
$$
}
\vsk>
\nt
In this special situation one can describe explicitly the set
of all points $h\in\ZZ$ where the operator $R_{\al\be}(h)$
is non-invertible. It is the main result of this section, cf.\ [DO].

\proclaim{Theorem 2.3}
The operator $R_{\al\be}(h)$ is not invertible at $h\in\ZZ$ if and only if
$$
\al_i\!+\!\be_N\ts,\al_{i+1}\!+\!\be_{N-1}\lcd\al_N\!+\!\be_i
\ts<\ts h+i\ts\le\ts
\al_1\!+\!\be_i\ts,\al_2\!+\!\be_{i-1}\lcd\al_i\!+\!\be_1
$$
for at least one index $i\in\{1\lcd N\}$.
\endproclaim

\demo{Proof}
By applying [\ts N,\,Proposition 2.12\ts] to the usual Young diagram $\be$ and
by employing the equality $b_1=0$ we can bring the product \(2.3) to the form
$$
\prod_k^\leftarrow\
\Bigl(1-\frac{P_{k,m+1}+\ldots+P_{k,m+n}}{h+a_k-z}\Bigr)\,\cdot\,
Y_\al\ox Y_\be
\nopagebreak
\Tag{2.4}
$$
where the factors corresponding to $k=1\lcd m$ are arranged from right to left.

\nt
Let us use the following well-known property of the element
$Y_\al\in\EndCN^{\ox m}$:
$$
\bigl(\ts P_{k,k+1}+\ldots+P_{km}\ts\bigr)\cdot Y_\al=-\ts a_k\,Y_\al\,,
\Tag{2.5}
$$
see [C1,\ts Theorem 4] and [\ts J,\,Section 4\ts].
For each $k=1\lcd m$ \text{consider the element}
$$
X_k=P_{k,k+1}+\ldots+P_{km}+P_{k,m+1}+\ldots+P_{k,m+n}
$$
of the algebra $\EndCN^{\ox(m+n)}$.
The sum $P_{k,k+1}+\ldots+P_{km}$ in $\EndCN^{\ox(m+n)}$
commutes with each of the elements $X_1\lcd X_{k-1}$.
Therefore by applying \(2.5) consecutively to $k=1\lcd m$
we can rewrite the product \(2.4) as
$$
\prod_k\ \Bigl(1-\frac{X_k+a_k}{h+a_k-z}\Bigr)\,\cdot\,Y_\al\ox Y_\be\,=\,
\prod_k\ \frac{h-X_k-z}{h+a_k-z}\,\cdot\,Y_\al\ox Y_\be\,.
\Tag{2.6}
$$
The elements $X_1\lcd X_m$ pairwise commute, so ordering of
the factors in \(2.6) corresponding to the indices $k$ does not matter.

Now for each $l=1\lcd n$ introduce the element
of the algebra $\EndCN^{\ox(m+n)}$
$$
X_{m+l}=P_{m+l,m+l+1}+\ldots+P_{m+l,m+n}\,.
$$
All the elements $X_1\lcd X_{m+n}$ also paiwise commute.
Let us use the properties of the element
$Y_\be\in\EndCN^{\ox n}$ similar to \(2.5)\ts:
$$
\bigl(\ts P_{1l}+\ldots+P_{l-1,l}\ts\bigr)\cdot Y_\be=b_l\,Y_\be\,;
\qquad
l=1\lcd n\,.
$$
Any symmetric polynomial in the sums $P_{12}+\ldots+P_{1l}$ with $l=1\lcd n$
belongs to the image in $\EndCN^{\ox n}$
of the centre of the group ring $\CC\cdot S_n$. Therefore
$$
\gather
\prod_l\,
\bigl(\ts h-P_{l,l+1}-\ldots-P_{ln}-z\ts\bigr)\,\cdot\,Y_\be\ =\
\\
\prod_l\,
\bigl(\ts h-P_{1l}-\ldots-P_{l-1,l}-z\ts\bigr)\,\cdot\,Y_\be\ =\
\prod_l\,
(\ts z-b_l-h\ts)\,\cdot\,Y_\be
\endgather
$$
in $\EndCN^{\ox n}$.
Hence we can rewrite the right hand side of the equiality \(2.6) as
$$
\prod_k\ \frac{h-X_k-z}{h+a_k-z}
\,\cdot\,
\prod_l\ \frac{h-X_{m+l}-z}{h-b_l-z}
\,\cdot\,
Y_\al\ox Y_\be\,.
\Tag{2.7}
$$
Note that for any $h\ts,z\in\CC$ the operator $(h-X_1-z)\ldots(h-X_{m+n}-z)$
is the image in $\EndCN^{\ox(m+n)}$ of a central
element in the group ring $\CC\cdot S_{m+n}$.
The eigenvalue of that central element in the irreducible {\=$S_{m+n}$}-module
corresponding to the Young diagram with the contents $c_1\lcd c_{m+n}$
equals $(h-c_1-z)\ldots(h-c_{m+n}-z)$.

The image of the operator $Y_\al\ox Y_\be$ in $(\CC^N)^{\ox(m+n)}$
is the tensor product of the irreducible {\=$\glN$}-modules of highest
weights $(\al_1\lcd\al_N)$ and $(\be_1\lcd\be_N)$.
Let $\ga=(\ga_1\lcd\ga_N)$ run through the set
$\Gal$ of highest weights
of all irreducible {\=$\glN$}-modules appearing in this tensor product.
We identify $\ga$ with corresponding usual Young diagram.
The expression \(2.7) for the product \(2.3) shows that the operator
$R_{\al\be}(h)$ is not invertible if and only if $h$ is a content for
at least one but not for all diagrams $\ga\in\Gal$. That is if and only if
$h\in\ZZ$ and
$$
\min_\ga\,\ga_i\,<h+i\,\le\,\max_\ga\,\ga_i
$$
for some $i\in\{1\lcd N\}$.
Let the index $i$ be fixed. Theorem 2.3 will follow from
$$
\align
\min_\ga\,\ga_i\,&=\,\max\,
\bigl(\ts
\al_i\!+\!\be_N\ts,\al_{i+1}\!+\!\be_{N-1}\lcd\al_N\!+\!\be_i
\ts\bigr)\,,
\Tag{2.8}
\\
\max_\ga\,\ga_i\,&=\,\min\,
\bigl(\ts
\al_1\!+\!\be_i\ts,\al_2\!+\!\be_{i-1}\lcd\al_i\!+\!\be_1
\ts\bigr)\,.
\Tag{2.9}
\endalign
$$

First let us check that the numbers $\ga_i$ with the fixed index $i$ attain
the values at the right hand sides of \(2.8) and \(2.9). We will use
the following well-known fact [\ts B,\ts Exercice\,VIII.9.14\ts]\,:
the unique dominant weight in the {\=$S_N$}-orbit
of the weight $(\al_1+\be_N\ts,\al_2+\be_{N-1}\lcd\al_N+\be_1)$
belongs to $\Gal$. Apply this fact to
$(\al_i\lcd\al_N)$ and $(\be_i\lcd\be_N)$ instead of
$(\al_1\lcd\al_N)$ and $(\be_1\lcd\be_N)$.
Then with the help of Littlewood-Richardson rule [\ts M\ts,\,Theorem\ I.9.2\ts]
one obtains that the set $\Gal$ contains the weight $\ga$ such that
$\ga_j=\al_j+\be_j$ for $j<i$ while
$$
\ga_i=
\max\,\bigl(\ts\al_i\!+\!\be_N\ts\lcd\al_N\!+\!\be_i\ts\bigr)\,,
\,\ \ldots\ \ts,\,
\ga_N=
\max\,\bigl(\ts\al_i\!+\!\be_N\ts\lcd\al_N\!+\!\be_i\ts\bigr)\,.
\nopagebreak
$$
So the value at the right hand side of \(2.8) is attained by $\ga_i$.
Similarly, by applying that well-known fact to
$(\al_1\lcd\al_i)$ and $(\be_1\lcd\be_i)$
one gets $\ga\in\Gal$ where
$$
\ga_1=
\max\,\bigl(\ts\al_1\!+\!\be_i\ts\lcd\al_i\!+\!\be_1\ts\bigr)\,,
\,\ \ldots\ \ts,\,
\ga_i=
\min\,\bigl(\ts\al_1\!+\!\be_i\ts\lcd\al_i\!+\!\be_1\ts\bigr)
\nopagebreak
$$
and $\ga_j=\al_j+\be_j$ for $j>i$.
\text{So $\!\ga_i\!$ attains the value at the right hand side of \(2.9).\!}

It now remains to prove for the fixed indices $i\ts,j\in\{1\lcd N\}$
the inequalities
$$
\alignat2
&\ga_i\,\ge\,\al_{N-j+i}+\be_j
&&\qquad\text{if}\qquad
i\le j\,,
\Tag{2.88}
\\
&\ga_i\,\le\,\al_{i-j+1}+\be_j
&&\qquad\text{if}\qquad
i\ge j\,.
\Tag{2.99}
\endalignat
$$
These inequalities follow easily from the Littlewood-Richardson rule.
Recall that the diagram $\ga$ occurs in $\Gal$ if there exists a semi-standard
tableau $\ka$ of shape $\ga/\be$ such that the values $1\lcd N$ are
taken by $\ka$ respectively $\al_1\lcd\al_N$ times and
the tableau $\ka$ satisfies the lattice property [\ts M\ts,\,Section\ I.9\ts].
Fix such a tableau $\ka$ and for $j\ts,k=1\lcd N$ denote by $\al_{jk}$ the
number of times $\ka$ takes the value~$k$ in the row $j$.
Then $\al_{jk}=0$ for $j<k$ since $\ka$ is semi-standard and has the
lattice property. Let us check inequality \(2.88).
Write $l=N-j+i$ for short. If $i\le j$
$$
\al_l=
\al_{ll}+\ldots+\al_{Nl}\le
\al_{l-1,l-1}+\ldots+\al_{N-1,l-1}\le
\,\ldots\,\le
\al_{ii}+\ldots+\al_{ji}
\nopagebreak
$$
because of the lattice property. Further, here we have
$
\al_{ii}+\ldots+\al_{ji}\le\ga_i-\be_j
$
since $\ka$ is semi-standard. Thus we obtain \(2.88).

Now suppose that $i\ge j$. Let us write $l=i-j+1$ for short. Then we have
$
\ga_i-\be_j\le\al_{ii}+\ldots+\al_{il}
$
because $\ka$ is semi-standard. Here by the lattice property
$$
\align
\al_{ii}+\ldots+\al_{il}\le
\al_{i-1,i-1}+\al_{i,i-1}+\ldots+\al_{il}&\le
\\
\al_{i-2,i-2}+\al_{i-1,i-2}+\al_{i,i-2}+\ldots+\al_{il}&\le
\,\ldots\,\le
\al_{ll}+\ldots+\al_{il}\le
\al_l\,.
\endalign
\nopagebreak
$$
Thus we obtain the inequality \(2.99) and complete the proof of Theorem 2.3.
\enddemos

%------------------------------------------------------------------------------

\bigskip
\Section{Cyclicity conditions}

\nt
In this section we will consider the tensor product of $n$ elementary
{\=$\YN$}-modules for any $n$. For each $s=1\lcd n$ fix a non-negative integer
$M\su$ and take a pair of
non-increasing sequences of integers $\la\su,\mu\su$ with
lengths $N+M\su,M\su$ respectively. Take a parameter $h\su\in\CC$.
Consider the elementary {\=$\YN$}-module $V\su=V_{\la\su,\ts\mu\su}(h\su)$.
Its Drinfeld polynomials $P_1\su(u)\lcd P_{N-1}\su(u)$ are given explicitly
by Proposition 1.6. Introduce the rational functions
$$
Q_k\su(u)=P_k\su(u)\big/P_k\su(u+1)\,;
\qquad
k=1\lcd N-1\,.
\Tag{3.333}
$$
Further, denote by $\X_k\su$ the collection of all numbers
of the form
$$
i-h\su-\la_{M\su+k,i}-1\,;
\qquad
i=1\lcd M\su+k
\Tag{3.111}
$$
where $\la_{M\su+k,i}$ is the {\=$(M\su+k,i)$}-entry of any
scheme $\La$ in $\S_{\la\su,\ts\mu\su}$ such that
the array $\La^\tm$ obtained from $\La$ by
decreasing this entry by $1$ is also in $\S_{\la\su,\ts\mu\su}$.

Let $\ze\su\in V\su$ be a singular vector,
it is determined up to scalar multiplier.
Consider {\=$\YN$}-module $V$ obtained from the tensor product
$V^{(1)}\ox\ldots\ox V^{(n)}$ via the comultiplication \(1.7).
Consider the vector $\ze=\ze^{(1)}\ox\ldots\ox\ze^{(n)}\in V$.
By Proposition 1.8 for any $k=1\lcd N-1$ the coproduct
$\Den\bigl(C_k(u)\bigr)$ is a sum of the elements in $\YN^{\ox n}\ts[[u^{-1}]]$
with the degrees in $\ZZ^n$ containing at least one positive component.
Therefore $C_k(u)\cdot\ze=0$ in the module $V$ for any index $k$. The next
proposition gives sufficient conditions for cyclicity of the vector $\ze$
under the action of the coefficients of the series $B_1(u)\lcd B_{N-1}(u)$
in the module $V$.

\proclaim{Proposition 3.1}
Suppose that $Q_k\su(x)\neq0$ for any $x\in\X_k\ru$ when $1\le k<n$
and $1\le r<s\le n$.
Then the vector $\ze$ in the {\=$\YN$}-module $V$ is cyclic.
\endproclaim

\demo{Proof}
Each of the vector spaces $V_{\la\su,\ts\mu\su}$ has a {\=$\ZZ$}-grading
as defined in the end of Section 1. Thus the space $V$ acquires grading
by the elements of the group $\ZZ^n$. We will equip the set $\ZZ^n$ with
lexicographical ordering: the element $(d_1\lcd d_n)$ precedes
$(d_1^{\,\prime}\lcd d_n^{\,\prime})$ if $d_s<d_s^{\,\prime}$ for some index
$s$ while $d_r=d_r^{\,\prime}$ for each $r<s$.

Take any vector $\xi=\xi^{(1)}\ox\ldots\ox\xi^{(n)}\in V$
where any tensor factor $\xi\su$ is an element of the Gelfand-Zetlin basis
in $V_{\la\su,\ts\mu\su}$. Fix any index $r\in\{1\lcd n\}$ and assume that
$\xi\su=\ze\su$ for every $s>r$.
Write $\xi\ru=\xi_\La$ for a certain scheme $\La$ in $\S_{\la\ru,\ts\mu\ru}$.
Then fix any indices $k\in\{1\lcd N-1\}$ and $i\in\{1\lcd M\ru+k\}$
such that the array $\La^\tm$ obtained from $\La$ by
decreasing the {\=$(M\ru+k,i)$}-entry \text{by $1$,}
is again in $\S_{\la\ru,\ts\mu\ru}$.
We have $x=i-h\ru-\la_{M\ru+k,i}-1\in\X_k\ru$.
Take
$$
\xi^{(1)}\ox\ldots\ox\xi^{(r-1)}
\ox\xi_{\La^\tm}\ox
\ze^{(r+1)}\ox\ldots\ox\ze^{(n)}
\in V\ts.
\Tag{3.0}
$$

Consider the rational function $B_k(u)\cdot\xi$ of $u$ valued in $V$.
Take the term of this function with the leading degree in $\ZZ^n$. It
is again a rational function of $u$ valued in $V$.
Let $\eta\in V$ be the first non-zero coefficient of the Laurent series
in $u-x$ of the latter function. Here and in what follows we consider
the Laurent expansions at $u=x$. We will show that the vector $\eta$
is a scalar multiple of \(3.0). This guarantees the cyclicity of the vector
$\ze=\ze^{(1)}\ox\ldots\ox\ze^{(n)}$ under the action of the coefficients
of the series $B_1(u)\lcd B_{N-1}(u)$ in the module $V$.

Observe that by Proposition 1.8 the coproduct $\De\bigl(B_k(u)\bigr)$
\text{is equal to the sum}
$$
A_k(u)\ox B_k(u)\,+\,B_k(u)\ox D_k(u)
\nopagebreak
$$
plus the terms of degrees in $\ZZ^2$ with a positive second component.
Therefore by our assumption on $\xi$ the vector
$B_k(u)\cdot\xi\in V^{(1)}\ox\ldots\ox V^{(n)}$ equals the sum
$$
\gather
\bigl(A_k(u)\cdot\xi^{(1)}\ox\ldots\ox\xi^{(r)}\bigr)\ox
\bigl(B_k(u)\cdot\ze^{(r+1)}\ox\ldots\ox\ze^{(n)}\bigr)\,+
\Tag{3.1}
\\
\bigl(B_k(u)\cdot\xi^{(1)}\ox\ldots\ox\xi^{(r)}\bigr)\ox
\bigl(D_k(u)\cdot\ze^{(r+1)}\ox\ldots\ox\ze^{(n)}\bigr)\ts,
\endgather
$$
where the actions of the algebra $\YN$ in
$V^{(1)}\ox\ldots\ox V^{(r)}$ and $V^{(r+1)}\ox\ldots\ox V^{(n)}$
are determined via the comultiplications
$\Der$ and $\De^{\hskip-1pt(n-r)}$ respectively.

%Let $\operatorname{J}$ be the left ideal in the algebra $\YN\ts[[u^{-1}]]$
%generated by elements with positive {\=$\ZZ$}-degrees.
By the second equality in \(1.7777777)
the series $D_k(u)$ is equal to the sum of
$$
T_{k+1,k+1}(u-k+1)\cdot
T_{k-1,k-1}(u-k+2)\ldots T_{22}(u-1)\,T_{11}(u)
$$
and certain elements from the ideal $\operatorname{J},$ see Section 1.
\text{\!The series $B_k(u)$ is the sum of\!}
$$
T_{k,k+1}(u-k+1)\cdot
T_{k-1,k-1}(u-k+2)\ldots T_{22}(u-1)\,T_{11}(u)
$$
and again of certain elements from $\operatorname{J}$. But the vector
$\ze^{(r+1)}\ox\ldots\ox\ze^{(n)}$ is an eigenvector for the action
in $V^{(r+1)}\ox\ldots\ox V^{(n)}$ of the product
$$
T_{k-1,k-1}(u-k+2)\ldots T_{22}(u-1)\,T_{11}(u)\,.
$$
Hence by dividing the vector \(3.1) by
the corresponding eigenvalue we get the sum
$$
\gather
\ \quad
\bigl(A_k(u)\cdot\xi^{(1)}\ox\ldots\ox\xi^{(r)}\bigr)\ox
\bigl(\ts T_{k,k+1}(u-k+1)\cdot\ze^{(r+1)}\ox\ldots\ox\ze^{(n)}\bigr)\,+
\Tag{3.2}
\\
\qquad
\bigl(B_k(u)\cdot\xi^{(1)}\ox\ldots\ox\xi^{(r)}\bigr)\ox
\bigl(\ts T_{k+1,k+1}(u-k+1)\cdot\ze^{(r+1)}\ox\ldots\ox\ze^{(n)}\bigr)\ts.
\endgather
$$
Note that by taking here the first non-zero Laurent coefficient
at $u=x$
of the component with the leading degree in $\ZZ^n$, we get a scalar multiple
of the same vector $\eta$ as in $B_k(u)\cdot\xi$.

Further, by Proposition 1.8 the components of
$\Der\bigl(A_k(u)\bigr)$ and $\Der\bigl(B_k(u)\bigr)$ with the
leading degrees in $\ZZ^n$ are respectively
$$
A_k(u)^{\ox(r-1)}\ox A_k(u)
\quad\text{and}\quad
A_k(u)^{\ox(r-1)}\ox B_k(u).
$$
But by Theorem 1.3 the tensor product
$
\bigl(A_k(u)\cdot\xi^{(1)}\bigr)\ox\ldots\ox\bigl(A_k(u)\cdot\xi^{(r-1)}\bigr)
$
equals $\xi^{(1)}\ox\ldots\ox\xi^{(r-1)}$
times a certain rational function of $u$ valued in $\CC$.
Divide \(3.2) by this rational function.
Now it suffices to show that by taking in
$$
\gather
\bigl(A_k(u)\cdot\xi_\La\bigr)\ox
\bigl(\ts T_{k,k+1}(u-k+1)\cdot\ze^{(r+1)}\ox\ldots\ox\ze^{(n)}\bigr)\,+
\Tag{3.3}
\\
\quad
\bigl(B_k(u)\cdot\xi_\La\bigr)\ox
\bigl(\ts T_{k+1,k+1}(u-k+1)\cdot\ze^{(r+1)}\ox\ldots\ox\ze^{(n)}\bigr)
\endgather
$$
the first non-zero Laurent coefficient at $u=x$,
we get scalar multiple of the vector
$$
\xi_{\La^\tm}\ox\ze^{(r+1)}\ox\ldots\ox\ze^{(n)}\ts.
$$

Let $R_{k+1}^{(r+1)}(u)\lcd R_{k+1}^{(n)}(u)$
be eigenvalues of $T_{k+1,k+1}(u)$ on $\ze^{(r+1)}\lcd\ze^{(n)}$
respectively. By \(1.6)
coproduct $\De^{\hskip-1pt(n-r)}\bigl(\ts T_{k+1,k+1}(u)\bigr)$
equals $T_{k+1,k+1}(u)^{\ox(n-r)}$ plus terms with the degrees in
$\ZZ^{n-r}$ containing at least one positive \text{component. So}
$$
T_{k+1,k+1}(u)
\cdot\ze^{(r+1)}\ox\ldots\ox\ze^{(n)}=
R_{k+1}^{(r+1)}(u)\ldots R_{k+1}^{(n)}(u)
\cdot\ze^{(r+1)}\ox\ldots\ox\ze^{(n)}.
$$
Determine the rational function $\rho_k(u)$ by \(1.007) for
\text{$M=M\ru,h=h\ru,\mu_i=\mu_i\ru$.}
By Theorem 1.3 the value of $\rho_k(u)\ts A_k(u)\cdot\xi_\La$
at $u=x$ is zero.
By Theorem 1.4 the value of
$\rho_k(u)\ts B_k(u)\cdot\xi_\La$ at $u=x$
is a non-zero scalar multiple \text{of vector $\xi_{\La^\tm}$.}
Divide \(3.3) by
$$
R_{k+1}^{(r+1)}(u-k+1)\ts\ldots\ts R_{k+1}^{(n)}(u-k+1)\,.
$$
It now remains to prove regularity at $u=x$ of the rational function %in $u$
$$
\frac
{\,T_{k,k+1}(u-k+1)\cdot\ze^{(r+1)}\ox\ldots\ox\ze^{(n)}}
{R_{k+1}^{(r+1)}(u-k+1)\ts\ldots\ts R_{k+1}^{(n)}(u-k+1)}
\,.
\Tag{3.4}
$$

Thanks to the equalities \(1.999999999) we get from \(3.333) that
for every $s=r+1\lcd n$
$$
Q_k\su(u)=R_{k+1}\su(u-k+1)\big/R_k\su(u-k+1)\,.
\Tag{3.5}
$$
On the other hand, the coproduct
$\De^{\hskip-1pt(n-r)}\bigl(\ts T_{k,k+1}(u)\bigr)$ has the form
$$
\sum_{s>r}\
T_{kk}(u)^{\ox(s-r-1)}\ox T_{k,k+1}(u)\ox T_{k+1,k+1}(u)^{\ox(n-s)}
\nopagebreak
$$
plus terms with the degrees in
$\ZZ^{n-r}$ containing at least one positive \text{component. So} due to \(3.5)
the vector \(3.4) is equal to the sum over $s=r+1\lcd n$ of the vectors
$$
\frac{\ze^{(r+1)}\ox\ldots\ox\ze^{(s-1)}}
{Q_k^{(r+1)}(u)\ts\ldots\ts Q_k^{(s-1)}(u)}
\ox
\frac{T_{k,k+1}(u-k+1)\cdot\ze\su}{R_{k+1}\su(u-k+1)}
\ox
\ze^{(s+1)}\ox\ldots\ox\ze^{(n)}\,.
$$
Here $Q_k^{(r+1)}(x)\lcd Q_k^{(s-1)}(x)\neq0$ by our assumption. Thus
to complete the proof it suffices to show that the rational function in $u$
$$
{T_{k,k+1}(u-k+1)\cdot\ze\su}\big/{R_{k+1}\su(u-k+1)}
\Tag{3.55}
$$
is also regular at $u=x$. Suppose this function is not identically zero.
Denote by $\varpi$ the first non-zero Laurent coefficient of this function
at $u=x$.
Let $v$ be a formal parameter.
The {\=$\YN$}-module $V\su$ is irreducible and its subspace of
maximal {\=$\ZZ$}-degree is spanned by vector $\ze\su$.
Therefore $T_{l+1,\ts l}(v)\cdot\varpi\neq0$ as a formal series in $v$
for at least one index $l\in\{1\lcd N-1\}$. But
$$
T_{l+1,\ts l}(v)\,T_{k,k+1}(u)\cdot\ze\su=
\frac{T_{kl}(v)\,T_{l+1,\ts k+1}(u)-T_{kl}(u)\,T_{l+1,\ts k+1}(v)}{u-v}
\cdot\ze\su
\Tag{3.6}
$$
due to the defining relations \(1.1). Hence $T_{l+1,\ts l}(v)\cdot\varpi=0$
for any index $l>k$. If $l<k$ then by applying \(1.1) to the right hand side
of the equality \(3.6), we get $T_{l+1,\ts l}(v)\cdot\varpi=0$ again.
Thus the series $T_{k+1,\ts k}(v)\cdot\varpi$ in $v$ is not identically zero.
On the other hand, by applying \(3.6) to $l=k$ we get the equality
$$
T_{k+1,k}(v)\cdot
\frac{T_{k,k+1}(u-k+1)}{R_{k+1}\su(u-k+1)}\cdot\ze\su=
\frac1{u-v-k+1}
\biggl(R_k\su(v)-\frac{R_{k+1}\su(v)}{Q_k\su(u)}\ts\biggr)\cdot\ze\su
$$
where the right hand has no pole at $u=x$ because
$Q_k\su(x)\neq0$ by our assumption.
Therefore the rational function \(3.55) is indeed regular at $u=x$.
\enddemos

\nt
The next proposition
matches Proposition 3.1 and is essentially equivalent to it.

\proclaim{Proposition 3.2}
Suppose that $Q_k\su(x)\neq0$ for any $x\in\X_k\ru$ when $1\le k<n$ and
$1\le s<r\le n$.
Then the vector $\ze$ in the {\=$\YN$}-module $V$ is cocyclic.
\endproclaim

\demo{Proof}
Let us consider the {\=$\YN$}-module dual to $V$. Due to Proposition 1.7
and to \(1.666666) it is equivalent to the {\=$\YN$}-module obtained from
the tensor product $V^{(1)}\ox\ldots\ox V^{(n)}$
by composing the comultiplication \(1.7) with the transposition
$(1\lcd n)\mapsto(n\lcd1)$ of the tensor factors.
Denote by $V^\prime$ the latter module.
The cocyclicity of the vector $\ze=\ze^{(1)}\ox\ldots\ox\ze^{(n)}$ in $V$
amounts to the cyclicity of the same vector in $V^\prime$.
Now Proposition 3.1 provides the required statement.
\enddemos

\nt
By combining Propositions 3.1 and 3.2 we immediately obtain
sufficient conditions for
irreducibility of the {\=$\YN$}-module $V=V^{(1)}\ox\ldots\ox V^{(n)}$.

\proclaim{Theorem 3.3}
Suppose that $Q_k\su(x)\neq0$ for any $x\in\X_k\ru$ whenever $1\le k<n$ and
$r\neq s$. Then the {\=$\YN$}-module $V$ is irreducible.
\endproclaim

\nt
In general, these conditions are not necessary for the irreducibility of $V$.
Still by using again Proposition 3.1 we can give a criterion for the
irreducibility of $V$ when each of the skew Young diagrams
$\la^{(1)}/\mu^{(1)}\lcd\la^{(n)}/\mu^{(n)}$ has the simplest shape.

First let us make a general remark.
For $s=1\lcd n$ take the {\=$\YN$}-module $V_{\la\su/\mu\su}(h)$
as defined in \text{Section 2.}
Due to Proposition 2.1 the {\=$\YN$}-module $V\su$ can be
obtained from it by pulling back through
an automorphism of the form $\of$. Let us fix this realization of $V\su$.
Note that by the definition \(1.6) we have the equalities
$$
\De\o\of=(\ts\of\ox\id\ts)\o\De=(\ts\id\ox\of\ts)\o\De\,.
$$
Therefore for any $r<s$ the element
$$
R_{\la\ru/\mu\ru,\ts\la\su/\mu\su}\bigl(\ts h\ru-h\su)
\in\End\bigl(\,\im Y_{\la\ru/\mu\ru}\ox\im Y_{\la\su/\mu\su}\bigr)
\Tag{4.1}
$$
is an intertwining operator between the {\=$\YN$}-modules obtained from
the tensor product $V\ru\ox V\su$ via the comultiplications $\Dep$ and $\De$
respectively.

\proclaim{Theorem 3.4}
Suppose that each of the skew
diagrams $\la^{(1)}/\mu^{(1)}\lcd\la^{(n)}/\mu^{(n)}$
has rectangular shape.
Then the {\=$\YN$}-module $V$ is irreducible if and only if the operator
\(4.1)
is invertible whenever $1\le r<s\le n$.
\endproclaim

\demo{Proof}
If at least one of the operators \(4.1) with $r<s$ is not invertible
then the {\=$\YN$}-module $V$ is reducible thanks to the general remark
made above. Suppose that each of the operators \(4.1) with $r<s$ is invertible.
We will show that then the {\=$\YN$}-module $V$ is irreducuble.
Without any loss of generality we can assume that
$$
V\su=V_{\la\su,\ts\mu\su}\bigl(h\su\bigr)=V_{\la\su/\mu\su}\bigl(h\su\bigr)
$$
for each $s=1\lcd n$.
Recall that if a skew diagram $\al$ is obtained by adding
the same number $c$ to every contents of a skew diagram $\be$ then
$V_\al(h)=V_\be(h+c)$. So we can further assume that for each $s=1\lcd n$
$$
\la\su=\bigl(k\su\lcd k\su,0\lcd 0\bigr)
\quad\text{and}\quad
\mu\su=\varnothing
\Tag{4.4}
$$
where the positive integer $k\su$ appears $l\su\le N$ times.
If $l\su$ equals $0$ or $N$ then the {\=$\YN$}-module $V\su$ is
one-dimensional. We will assume that $0<l\su<N$.

If for some $k\in\{1\lcd N-1\}$
the number \(3.111) appears in the set $\X_k\su\!$ then
$$
\max\ts(\ts0\ts,l\su-N+k\ts)<i\le\min\ts(\ts k,l\su)
$$
and in this case the integer $\la_{M\su+k,i}=\la_{ki}$
can vary from $1$ up to $k\su$. Therefore
$$
\X_k\su=\{\ts h-h\su-1\,\ts|\ts
\max(\ts0\ts,l\su-N+k\ts)-k\su\!<h<\min(\ts k,l\su)\ts,\,
h\in\ZZ\,\}\,.
$$
By Proposition 1.6 we get $P_k\su(u)=Q_k\su(u)=1$ for any $k\neq l\su$.
If $k=l\su$ then
$$
P_k\su(u)=(\ts u+h\su-l\su+1)\ldots(\ts u+h\su-l\su+k\su)
\nopagebreak
$$
and the rational function $Q_k\su(u)$ has the only one zero $u=l\su-h\su-1$.

Let us fix any indices $r<s$.
Then $Q_k\su(u)$ has a zero in $\X_k\ru$ only for $k=l\su$ and only when
we have the inequalities
$$
-\min\bigl(\ts l\su\!\ts,N-l\ru\bigr)-k\ru
\!<h\ru-h\su\!<
\min\bigl(\ts0\ts,l\ru-l\su\bigr)
\Tag{4.2}
$$
while $h\ru-h\su\in\ZZ$. By exchanging the triples
$(\ts h\ru,k\ru,l\ru)$ and $(\ts h\su,k\su,l\su)$ in \(4.2) we
obtain the inequalities
$$
\max\bigl(\ts0\ts,l\ru-l\su\bigr)
\!<h\ru-h\su\!<
\min\bigl(\ts l\ru\!\ts,N-l\su\bigr)+k\su
\nopagebreak
\Tag{4.3}
$$
where again $h\ru-h\su\in\ZZ$. Now observe that the inequalities
\(4.2) and \(4.3) exclude each other. On the other hand, the {\=$\YN$}-modules
$V\ru\ox V\su$ and $V\su\ox V\ru$ obtained via the comultiplication $\De$
are equivalent: composition of the exchange map
$V\su\ox V\ru\to V\ru\ox V\su$ with \(4.1) is invertible and commutes
with the action of $\YN$. So we can assume that
$Q_k\su(x)\neq0$ for any $x\in\X_k\ru$ and $k\in\{1\lcd n-1\}$.
Then the vector $\ze\in V$ is cyclic by Proposition 3.1.

The {\=$\YN$}-module $V^\prime$ introduced in the proof of
Corollary 3.2 is equivalent to $V$. The isomorphism $V^\prime\to V$
is given by composition of the operators \(4.1) with
$$
(r,s)=(1,2)\ts,(1,3)\ts,(2,3)\lcd\ldots\lcd(1,n)\lcd(n-1,n)\ts.
$$
This isomorphism preserves one-dimensional subspace in $V$ spanned
by vector $\ze$. Indeed, each tensor factor $\ze\su\in V\su$ of
$\ze=\ze^{(1)}\ox\ldots\ox\ze^{(n)}$
has the maximal degree with respect to {\=$\ZZ$}-grading by eigenvalues
of the action in $V\su$ of the element
$$
N\ts E_{11}+(N-1)\ts E_{22}+\ldots+E_{NN}
\ts\in\,\glN\subset\ts\YN.
$$
But the operator \(4.1) commutes with the action of Lie algebra $\glN$ in
\text{$V\ru\ox V\su.$\!}
So cyclicity of the vector $\ze$ in the module $V$
is equivalent to its cocyclicity.
\enddemos

\nt
When the diagrams $\la^{(1)}/\mu^{(1)}\lcd\la^{(n)}/\mu^{(n)}$ have rectangular
shapes, Theorem 2.3 explicitly describes for each $r<s$ the set of all points
$h\ru-h\su\in\ZZ$ where the operator \(4.1) is not invertible.
Under the assumptions \(4.4) the non-invertibility occurs
if and only if one of the next two pairs of inequalities holds:
$$
\gather
-\min\bigl(\ts l\su\!\ts,N-l\ru\bigr)\!-\!\ts k\ru
\!<h\ru-h\su\!<
\min\bigl(\ts0\ts,l\ru\!-l\su\bigr)\ts\!+\!\ts\min\bigl(\ts0\ts,k\su-k\ru\bigr)
\ts,
\\
\max\bigl(\ts0\ts,l\ru\!-l\su\bigr)\ts\!+\!\ts\max\bigl(\ts0\ts,k\su-k\ru\bigr)
\!<h\ru-h\su\!<
\min\bigl(\ts l\ru\!\ts,N-l\su\bigr)\!+\!\ts k\su
\ts.
\endgather
$$
Note that if $k\ru=k\su$ then these pairs coincide
with \(4.2) and \(4.3) respectively. Therefore if $k^{(1)}=\ldots=k^{(n)}$
then already the conditions of Theorem 3.3 are necessary and sufficient
for the irreducibility of the {\=$\YN$}-module $V$. In the case
$k^{(1)}=\ldots=k^{(n)}=1$ our Theorem 3.4 follows from
[\ts AK\ts,Theorem 4.1\ts]. In the other special case
$l^{(1)}=\ldots=l^{(n)}=1,$ Theorem 3.4 follows from
\text{[\ts Z\ts,Theorem 4.2\ts].}

%------------------------------------------------------------------------------

\section{Acknowledgements}

\nt
This work was done while the second author stayed at Mathematics
Department, Osaka University and the first author visited RIMS,
Kyoto University. We are very grateful to E.\,Date and T.\,Miwa
for discussions and hospitality at these institutes.
We are also grateful to G.\,Olshanski and A.\,Zelevinsky for helpful remarks.
The first author was supported by an EPSRC Advanced Research Fellowship.
The second author would like to thank Forschungsinstitut f\"ur Mathematik,
ETH Z\"urich which he visited while finishing this article.

\newpage

%------------------------------------------------------------------------------

\line{\bf Addendum\hfill}
\medskip\smallskip\nt
The purpose of this additional section is to derive from Proposition 3.2 
the theorem stated below. In different but equivalent forms was
conjectured in [C1] and proved in [AK]. Recall that for any
fundamental {\=$\glN$}-weight $\nu=(1\lcd 1,0\lcd 0)$ and for any $h\in\CC$
the elementary {\=$\YN$}-module
$V_{\nu,\raise1pt\hbox{$\hskip-1pt{}_\varnothing\hskip-2pt$}}(h)$ 
is also called \emph{fundamental}.
The theorem is particularly important.
Due to [D2] it implies that
every irreducible finite-dimensional {\=$\YN$}-module
is equivalent to a submodule of a tensor product of fundamental modules, 
pulled back through an automorphism of %the algebra 
$\YN$ of the form $\of$. 
Therefore our Proposition 3.2 also implies this equivalence property.

%Yet another proof of the theorem
%was given in [KN] by using [KNV].

To state the theorem, for $s=1\lcd n$ let 
$V\su=V_{\nu\su,\raise1pt\hbox{$\hskip-1pt{}_\varnothing\hskip-2pt$}}(h\su)$
where $\nu\su$ is any fundamental {\=$\glN$}-weight and $h\su\in\CC$.
Let $V$ and $V^\prime$ be the {\=$\YN$}-modules obtained
from the tensor product $V^{(1)}\ox\ldots\ox V^{(n)}$
by using the comultiplications $\De$ and $\Dep$ respectively.
For $s<r$ let $Z$ and $Z^\prime$ 
be the {\=$\YN$}-modules obtained from %the tensor product 
$V\su\ox V\ru$ by using $\De$ and $\Dep$
respectively. Define an intertwining operator
$R^{(sr)}:Z'\to Z$ as in our Introduction. Let 
$R:V^\prime\to V$ be the composition of operators acting 
non-trivially only on the $s$th and $r$th tensor factors, 
as $R^{(sr)}$ where
$
(s,r)=(1,2)\ts,(1,3)\ts,(2,3)\lcd\ldots\lcd(1,n)\lcd(n-1,n)\ts.
$

\proclaim{Theorem}
Suppose that for all $s<r$ the difference
$h\su-h\ru$ is not a positive integer.
Then %the intertwining operator 
the image of $R$ is an irreducible submodule of the 
{\=$\YN$}-module~$V$.
\endproclaim

\demo{Proof}
Take any two indices $s$ and $r$ such that $1\le s<r\le n$. 
Let $a$ be the multiplicity %numbers of appearance 
of the entry $1$ in %the sequences 
$\nu\su$.
By Proposition 1.6 for $k=1\lcd N-1$ the Drinfeld polynomial
$P_k\su(u)$ of the {\=$\YN$}-module $V\su$ is $u+h\su+1-a$ if $k=a$,
and is $1$ if $k\neq a$.
So the rational function 
$Q_k\su(u)$ has a zero only at $u=a-1-h\su$
and only if $k=a$, see the definition \(3.333). 

Now consider the set $\X_{\ts a}\ru$, see \(3.111). Take any
$\La$ in $\S_{\ts\nu\ru,\raise1pt\hbox{$\hskip-1pt{}_\varnothing\hskip-2pt$}}$.
If the array $\La^\tm$ obtained by
decreasing an entry of $\La$ by $1$ is also in
$\S_{\ts\nu\ru,\raise1pt\hbox{$\hskip-1pt{}_\varnothing\hskip-2pt$}}$ 
then this entry is $1$. So the set
$\X_{\ts a}\ru$ is contained in the set of numbers
$i-2-h\ru$ with $i=1\lcd a$. None of these numbers equals
$a-1-h^{(s)}$ because $h\su-h\ru\neq1\lcd a$.

By our Proposition 3.2 the vector
$\ze=\ze^{(1)}\ox\ldots\ox\ze^{(n)}$ in the {\=$\YN$}-module $V$ is cocyclic.
Here %the tensor factors
$\ze^{(1)}\lcd\ze^{(n)}$ denote singular vectors
in the {\=$\YN$}-modules $V^{(1)}\lcd V^{(n)}$ respectively.
They are determined %uniquely 
up to scalar multipliliers. Also 
the vector $\ze$ is cyclic in the {\=$\YN$}-module $V^{\prime}$,
see the proof of Proposition 3.2.

For $1\le s<r\le n$ the operator $R^{(sr)}$ on 
the vector space $V\su\ox V\ru$
preserves the one-dimensional subspace spanned
by vector $\ze\su\ox\ze\ru$, see the end of proof of Theorem 3.4.
%for a similar argument. 
To prove the above theorem, it now suffices to show that the vector
$R^{(sr)}\cdot\ze\su\ox\ze\ru\neq0$. 
Indeed, then $R\cdot\ze\neq0$ and %the vector 
$\ze$ is a multiple of
$R\cdot\ze$. 

The inequality  $R^{(sr)}\cdot\ze\su\ox\ze\ru\neq0$ follows from the cyclicity of %the vector
$\ze\su\ox\ze\ru$ in the $\YN$-module $Z^{\prime}$,
which a particular case of the cyclicity of $\ze$ in $V^{\prime}$ for $n=2\,$.
If $R^{(sr)}\cdot\ze\su\ox\ze\ru=0$ then the intertwining operator $R^{(sr)}:Z'\to Z$ 
would be zero, thus contradicting its definition.
\enddemos

In the course of the proof
we showed that under the conditions of the theorem, the vector $\ze$
belongs to the image of $R$. Hence the vector $\ze$ is singular for
this {\=$\YN$}-module. Thus for any $k=1\lcd N-1$ 
the Drinfeld polynomial $P_k(u)$ of the image of $R$ equals the product 
over $s=1\lcd n$ of the polynomials $P_k\su(u)$. Namely,
$P_k(u)$ equals the product of the factors
$u+h\su+1-k$ over all indices 
$s$ such that the multiplicity of the entry $1$ in $\nu\su$ is $k$. 

%------------------------------------------------------------------------------

\newpage

\section{References}

\itemitem{[AK]}
{T. Akasaka and M. Kashiwara},
\emph{Finite-dimensional representations of quantum affine algebras},
{Publ. Res. Inst. Math. Sci.}
{\bf 33}
(1997),
839--867.

\itemitem{[B]}
{N. Bourbaki},
\emph{Groupes at Alg\`ebres de Lie},
Ch.\,VII-VIII\ts, Hermann, Paris, 1975.

\itemitem{[C1]}
{I. Cherednik},
\emph{Special bases of irreducible representations of a degenerate
affine Hecke algebra},
{Funct.\ Analysis Appl.}
{\bf 20}
(1986),
76--78.

\itemitem{[C2]}
{I. Cherednik},
\emph{A new interpretation of Gelfand-Zetlin bases},
{Duke Math.\ J.}
{\bf 54}
(1987),
563--577.

\itemitem{[DO]}
{E. Date and M. Okado},
\emph{Calculation of excitation spectra of the spin model related with the
vector representation of the quantized affine algebra of type $A_n^{(1)}$},
{Internat.\ J.\ Modern Phys.}
{\bf A9}
(1994),
399--417.

\itemitem{[D1]}
{V. Drinfeld},
\emph{Hopf algebras and the quantum Yang-Baxter equation},
{Dokl.\ Math.}
{\bf 32}
(1985),
254--258.

\itemitem{[D2]}
{V. Drinfeld},
\emph{A new realization of Yangians and quantized affine algebras},
{Dokl.\ Math.}
{\bf 36}
(1988),
212--216.

\itemitem{[GZ]}
{I. Gelfand and M. Zetlin},
\emph{Finite-dimensional representations of the unimodular group},
{Dokl. Akad. Nauk SSSR}
{\bf 71}
(1950),
825--828.

\itemitem{[JK]}
{G. James and A. Kerber},
\emph{The Representation Theory of the Symmetric Group},
Addison-Wesley, Reading MA, 1981.

\itemitem{[J]}
{A. Jucys},
\emph{Symmetric polynomials and the centre of the symmetric group ring},
{Rep.\ Math.\ Phys.}
{\bf 5}
(1974),
107--112.

\itemitem{[KKN]}
{A. Kirillov, A. Kuniba and T. Nakanishi},
\emph{Skew Young diagram method
in spectral decomposition of integrable lattice models},
{Commun.\ Math.\ Phys.}
{\bf 185}
(1997),
441--465.

\itemitem{[M]}
{I. Macdonald},
\emph{Symmetric Functions and Hall Polynomials},
Clarendon Press, Oxford, 1979.

\itemitem{[MNO]}
{A. Molev, M. Nazarov and G. Olshanski},
\emph{Yangians and classical Lie algebras},
{Russian Math.\ Surveys}
{\bf 51}
(1996),
205--282.

\itemitem{[N]}
{M. Nazarov},
\emph{Yangians and Capelli identities},
{Amer.\ Math.\ Soc.\ Transl.}
{\bf 181}
(1997),
139--164.

\itemitem{[NT]}
{M. Nazarov and V. Tarasov},
\emph{Representations of Yangians with Gelfand-Zetlin bases},
{J.\ Reine Angew.\ Math.}
{\bf 496}
(1998),
181--212.

\itemitem{[O]}
{G. Olshanski},
\emph{Representations of infinite-dimensional classical groups,
limits of enveloping algebras, and Yangians},
{Adv. Soviet Math.}
{\bf 2}
(1991),
1--66.

\itemitem{[TU]}
{K. Takemura and D. Uglov},
\emph{The orthogonal eigenbasis and norms of eigenvectors in the spin
Calogero-Sutherland model},
{J.\ Phys}
{\bf A30}
(1997),
3685--3717.

\itemitem{[T]}
{V. Tarasov},
\emph{Irreducible monodromy matrices for the R-matrix of the XXZ-model
and local lattice quantum Hamiltonians},
{Theor. Math. Phys.}
{\bf 63}
(1985),
440--454.

\itemitem{[Z]}
{A. Zelevinsky},
\emph{Induced representations of reductive {\=$\frak{p}$}-adic groups
II. On irreducible representations of $GL(n)$},
{Ann. Scient. Ec. Norm. Sup.}
{\bf 13}
(1980),
165--210.

%------------------------------------------------------------------------------

\bye